\documentclass[a4paper,times]{article}
\usepackage[left=2.5cm,right=2.5cm,top=2.5cm,bottom=2.5cm]{geometry}
\usepackage{amsmath,amssymb,amscd,amsxtra,amsthm,amsfonts}
\usepackage{epsf,graphicx,epsfig,xcolor,latexsym,multirow}
\usepackage{enumerate}
\usepackage{mathrsfs}
\usepackage{enumitem}
\usepackage{framed}
\usepackage{diagbox}
\usepackage[caption=false]{subfig}
\usepackage{algorithmic}

\newtheorem{theorem}{Theorem}
\newtheorem{lemma}[theorem]{Lemma}
\newtheorem{remark}[theorem]{Remark}
\newtheorem{algorithm}[theorem]{Algorithm}
\newtheorem{example}[subsection]{Example}

\usepackage{amsopn}

\providecommand{\Div}{\operatorname{div}}          
          
\providecommand{\curl}{\operatorname{{\bf curl}}}  

\providecommand*{\Dist}[2]{\operatorname{dist}({#1};{#2})}   
\providecommand*{\Dist}[2]{\Dist{#1}{#2}}

\newcommand{\Vb}{{\mathbf{b}}}

\newcommand{\Ve}{{\mathbf{e}}}

\newcommand{\Vr}{{\mathbf{r}}}

\newcommand{\Vx}{{\mathbf{x}}}

\newcommand{\Bc}{{\boldsymbol{c}}}
\newcommand{\Bd}{{\boldsymbol{d}}}

\newcommand{\Bf}{{\boldsymbol{f}}}

\newcommand{\Bn}{{\boldsymbol{n}}}

\newcommand{\Bu}{{\boldsymbol{u}}}
\newcommand{\Bv}{{\boldsymbol{v}}}
\newcommand{\Bw}{{\boldsymbol{w}}}
\newcommand{\Bx}{{\boldsymbol{x}}}

\newcommand{\BA}{{\boldsymbol{A}}}
\newcommand{\BB}{{\boldsymbol{B}}}
\newcommand{\BC}{{\boldsymbol{C}}}
\newcommand{\BD}{{\boldsymbol{D}}}
\newcommand{\BE}{{\boldsymbol{E}}}

\newcommand{\BH}{{\boldsymbol{H}}}

\newcommand{\BJ}{{\boldsymbol{J}}}

\newcommand{\BL}{{\boldsymbol{L}}}

\newcommand{\BP}{{\boldsymbol{P}}}

\newcommand{\BU}{{\boldsymbol{U}}}
\newcommand{\BV}{{\boldsymbol{V}}}
\newcommand{\BW}{{\boldsymbol{W}}}

\newcommand{\Cf}{\mathcal{F}}

\newcommand{\Ct}{\mathcal{T}}

\newcommand{\bbA}{\mathbb{A}}
\newcommand{\bbB}{\mathbb{B}}
\newcommand{\bbC}{\mathbb{C}}
\newcommand{\bbD}{\mathbb{D}}
\newcommand{\bbE}{\mathbb{E}}
\newcommand{\bbF}{\mathbb{F}}
\newcommand{\bbG}{\mathbb{G}}
\newcommand{\bbH}{\mathbb{H}}
\newcommand{\bbI}{\mathbb{I}}
\newcommand{\bbJ}{\mathbb{J}}
\newcommand{\bbK}{\mathbb{K}}
\newcommand{\bbL}{\mathbb{L}}
\newcommand{\bbM}{\mathbb{M}}

\newcommand{\bbP}{\mathbb{P}}

\newcommand{\bbR}{\mathbb{R}}
\newcommand{\bbS}{\mathbb{S}}

\newcommand{\bbU}{\mathbb{U}}

\newcommand{\bbX}{\mathbb{X}}

\newcommand*{\SP}[2]{\left<{#1},{#2}\right>} 

\newcommand*{\N}[1]{\left\|{#1}\right\|}     
\newcommand*{\SN}[1]{\left|{#1}\right|}      
\renewcommand*{\SP}[2]{\left({#1},{#2}\right)} 

\newcommand{\mesh}{\Ct_h}

\newcommand*{\Lp}[2][\defaultdomain]{L^{#2}({#1})}
\newcommand*{\Lpv}[2][\defaultdomain]{\BL^{#2}({#1})}
\newcommand*{\NLp}[3][\defaultdomain]{\N{#2}_{\Lp[#1]{#3}}}
\newcommand*{\NLpv}[3][\defaultdomain]{\N{#2}_{\Lpv[#1]{#3}}}

\newcommand*{\Ltwo}[1][\defaultdomain]{\Lp[#1]{2}}
\newcommand*{\Ltwov}[1][\defaultdomain]{\Lpv[#1]{2}}

\newcommand*{\NLtwo}[2][\defaultdomain]{\NLp[#1]{#2}{2}}
\newcommand*{\NLtwov}[2][\defaultdomain]{\NLpv[#1]{#2}{2}}

\newcommand*{\Hm}[2][\defaultdomain]{H^{#2}({#1})}
\newcommand*{\Hmv}[2][\defaultdomain]{\BH^{#2}({#1})}
\newcommand*{\bHm}[3][\defaultdomain]{H_{#3}^{#2}({#1})}

\newcommand*{\bHmv}[3][\defaultdomain]{\BH_{#3}^{#2}({#1})}

\newcommand*{\Hone}[1][\defaultdomain]{\Hm[#1]{1}}
\newcommand*{\Honev}[1][\defaultdomain]{\Hmv[#1]{1}}

\newcommand*{\zbHone}[1][\defaultdomain]{\bHm[#1]{1}{0}}
\newcommand*{\zbHonev}[1][\defaultdomain]{\bHmv[#1]{1}{0}}

\newcommand*{\NHonev}[2][\defaultdomain]{{\N{#2}}_{\Honev[{#1}]}}
\newcommand*{\SNHone}[2][\defaultdomain]{{\SN{#2}}_{\Hone[{#1}]}}
\newcommand*{\SNHonev}[2][\defaultdomain]{{\SN{#2}}_{\Honev[{#1}]}}

\newcommand*{\Hdiv}[1][\defaultdomain]{\boldsymbol{H}(\Div,{#1})}
\newcommand*{\bHdiv}[2][\defaultdomain]{\boldsymbol{H}_{#2}(\Div,{#1})}
\newcommand*{\zbHdiv}[1][\defaultdomain]{\bHdiv[#1]{0}}
\newcommand*{\kHdiv}[1][\defaultdomain]{\boldsymbol{H}(\Div0,{#1})}

\newcommand*{\Hcurl}[1][\defaultdomain]{\boldsymbol{H}(\curl,{#1})}
\newcommand*{\bHcurl}[2][\defaultdomain]{\boldsymbol{H}_{#2}(\curl,{#1})}
\newcommand*{\zbHcurl}[1][\defaultdomain]{\bHcurl[#1]{0}}
\newcommand*{\kHcurl}[1][\defaultdomain]{\boldsymbol{H}(\curl0,{#1})}

\newcommand*{\NHcurl}[2][\defaultdomain]{\N{#2}_{\Hcurl[#1]}}

\newcommand*{\jump}[1]{\hspace{0.5mm}\left[\hspace{-1.5mm}\left[\hspace{+0.5mm}{#1}\hspace{+0.5mm}\right]\hspace{-1.5mm}\right]\hspace{0.5mm}}     

\newcommand{\be}{\begin{eqnarray}}
\newcommand{\ee}{\end{eqnarray}}

\newcommand{\ben}{\begin{eqnarray*}}
\newcommand{\een}{\end{eqnarray*}}

\newcommand*{\Nsemi}[1]{\SN{#1}_{1,h}}
\newcommand*{\Ndg}[1]{\N{#1}_{1,h}}

\newcommand\keywordsname{Key words}
\newcommand\AMSname{AMS subject classifications}

\newenvironment{@abssec}[1]
{\if@twocolumn
\section*{#1}%
\else
\vspace{.05in}\footnotesize
\parindent .2in
{\upshape\bfseries #1. }\ignorespaces
\fi}

{\if@twocolumn\else\par\vspace{.1in}\fi}

\newenvironment{keywords}{\begin{@abssec}{\keywordsname}}{\end{@abssec}}

\title{A constrained transport divergence-free finite element method for Incompressible MHD equations}
\author{Lingxiao Li
\thanks{Institute of Applied Physics and Computational Mathematics, Beijing, 100094, China.}
\and Donghang Zhang
\thanks{LSEC, NCMIS, Institute of Computational Mathematics and Scientific/Engineering Computing,
Academy of Mathematics and Systems Science, Chinese Academy of Sciences, Beijing 100190, China;
School of Mathematical Science, University of Chinese Academy of Sciences, Beijing 100049, China.}
\and Weiying Zheng
\thanks{LSEC, NCMIS, Institute of Computational Mathematics and Scientific/Engineering Computing,
Academy of Mathematics and Systems Science, Chinese Academy of Sciences, Beijing 100190, China;
School of Mathematical Science, University of Chinese Academy of Sciences, Beijing 100049, China.}}

\begin{document}
\date{}
\maketitle

\begin{abstract}
  In this paper we study finite element method for three-dimensional incompressible resistive magnetohydrodynamic equations,
  in which the velocity, the current density, and the magnetic induction are divergence-free. It is desirable that the discrete solutions should also satisfy divergence-free conditions exactly especially for the momentum equations. Inspired by constrained transport method,
  we devise a new stable mixed finite element method that can achieve the goal. We also prove the well-posedness of the discrete solutions.
  To solve the resulting linear algebraic equations, we propose a GMRES solver with an augmented Lagrangian block preconditioner.
  By numerical experiments, we verify the theoretical results and demonstrate the quasi-optimality of the discrete solver with respect to the number of degrees of freedom.
\end{abstract}


\section{Introduction}

The magnetohydrodynamic (MHD) equations use fluid theory to describe the interaction of charged particle under the influence of magnetic field. They have broad scientific and engineering applications, such as magnetic
fusion \cite{jardin2010}, astrophysics \cite{hawley1995,ruan2019} and liquid metals \cite{gerbeau2006, davidson2001}.
The MHD model is a typical multi-scale and multi-physics system. The flow of conducting fluid modifies electromagnetic fields, and conversely the electromagnetic fields modify the momentum of fluid through Lorentz force. The strong coupling between fluid and electromagnetic fields makes the design of effective numerical methods and scalable iterative
solvers very difficult (see e.g. \cite{ni2007II,phi16,li17,gre10,philip-chacon,chacon2008} and references therein).

In this paper, we study the stationary incompressible resistive MHD equations in a bounded,
simply-connected, and Lipschitz polyhedral domain $\Omega\subset\bbR^3$.
They comprise incompressible Navier-Stokes equations and stationary Maxwell's equations
\begin{subequations}\label{sys:MHD}
\begin{align}
\rho \Bu\cdot\nabla\Bu+\nabla{p}-\nu\Delta\Bu-\BJ\times\BB=\Bf\quad&\text{in}\;\;\Omega,\label{eq:NS}\\
\curl\BH=\BJ,\quad
\curl\BE=\mathbf{0}\quad&\text{in}\;\;\Omega,\label{eq:Faraday}\\
\Div\Bu=0,\quad\Div\BB=0, \quad \Div\BJ = 0\quad&\text{in}\;\;\Omega,\label{eq:Gauss}
\end{align}
\end{subequations}
where $\rho$ is the fluid density, $\Bu$ the fluid velocity, $p$ the hydrodynamic pressure, $\BE$ the electric feild, $\BH$ the magnetic field, $\BB$ the magnetic induction, $\BJ$ the electric current density, and $\Bf\in\Ltwov$ external force. The equations in (\ref{sys:MHD}) are complemented with the following constitutive equation and generalized Ohm's law
\begin{align}\label{eq:Ohm}
\BB=\mu\BH,\qquad\BJ=\sigma(\BE+\Bu\times\BB).
\end{align}
The physical parameters are, respectively, dynamic viscosity $\nu$, magnetic permeability $\mu$, and electric conductivity $\sigma$. For the well-posedness of (\ref{sys:MHD}) and (\ref{eq:Ohm}), we assume homogeneous Dirichlet boundary conditions
\begin{align}\label{eq:Bdry}
\Bu=\mathbf{0},\quad\BH\times\Bn=\mathbf{0}
\quad\text{on}\;\;\Gamma:=\partial\Omega.
\end{align}

For computational MHD, it is desirable to study discrete methods which respect the divergence-free constraint
for magnetic induction, namely $\Div\BB_h = 0$ \cite{brackbill1980,toth2000}.
There are already many important works on this topic in the literature.
For early works one can refer to the review paper \cite{toth2000} and the references therein,
where seven schemes are reviewed in detail, such as the 8-wave formulation,
the projection scheme, constrained transport (CT) methods and so on.
One important observation of the work is the close relation between
vector potential method (cf. e.g. \cite{salah2001}) and
CT methods (see \cite{evans1988}). In \cite{rossmanith2006}, Rossmanith proposed
an unstaggered and high-resolution CT method for MHD flows. The readers are also
referred to more recent papers \cite{felker2018, balsara2015, hu17, hip18} which deal with $\Div\BB_h=0$ exactly.
The second constraint to be satisfied is $\Div \Bu_h = 0$, that is,
mass conservation of fluid. Under certain extreme situations,
nonphysical phenomena may appear if discrete solutions are
not mass-conservative (see \cite{john2017} for comprehensive discussions).

The third constraint is $\Div \BJ_h = 0$, that is, charge conservation.
Assuming small magnetic Reynolds number, the full MHD equations can reduce to
a inductionless MHD equations, see \cite{ni2007II, ni2012III, zhang2014, li2019SIAM}. In this special case,
when the applied magnetic filed is constant, the discrete Lorentz force in the momentum equations,
which is a volume force, can only precisely conserve the total momentum when the current density is divergence-free.
The authors in \cite{ni2007II, ni2012III, zhang2014} suggested that only the numerical schemes,
which can conserve the total momentum in the discrete level, can obtain accurate result for MHD flows,
especially at large Hartmann numbers.
Based on the above considerations, we are motivated to develop a stable numerical scheme
which satisfies the three divergence-free conditions simultaneously in the \textit{momentum equation}, namely,
\begin{equation}\label{eq:maincon}
\Div \Bu_h = 0, \quad \Div \BJ_h = 0,\quad \Div \BB_h = 0.
\end{equation}
holds in the discrete scheme of momentum equation \eqref{eq:NS}.

In this paper, we propose a monolithic CT finite element method for \eqref{sys:MHD} such that the discrete
solutions in the momentum equations \eqref{eq:NS} satisfy \eqref{eq:maincon} exactly.
Compared with traditional CT methods, the new method treats magnetic field $\BH$ and magnetic vector potential $\BA$ as individual variables by means of edge finite element discretization \cite{ned86}, which doesn't need a staggered grid.
The discrete current density $\BJ_h:=\curl\BH_h$ and the discrete magnetic induction $\BB_h:=\curl\BA_h$ are divergence-free naturally.
As remarked in \cite{brackbill1980, evans1988,ni2007II,ni2012III}, numerical methods satisfying $\Div\BB_h = \Div\BJ_h = 0$
in the momentum equations can reduce nonphysical effects created by the discrete Lorentz force $\BJ_h\times\BB_h$ on the fluid movement.
To fulfill the first constraint in \eqref{eq:maincon},
we adopt precisely the same ideas as in \cite{cockburn2007,gre10} where
the velocity $\Bu$ and the pressure $p$ are discretized by $\Hdiv$-conforming
finite elements and fully discontinuous finite elements respectively.
The standard mixed finite element methods can ensure
$\|\Div\Bu_h\|_{L^2} = 0$ on the discrete level \cite{gre10}.

Now we introduce some but not at all complete references on finite element methods for incompressible MHD equations. More work can be found in the references of the cited papers. In \cite{salah2001},
Salah, Soulaimani, and Habashi used $(\Bu,p,\BB,\BA,\psi)$ as solution variables where $\psi$ is the scalar electric potential,
while adopted Galerkin-least-squares variational formulation and continuous finite element discretization.
In \cite{gerbeau2000}, Gerbeau developed a stabilized finite element method for the incompressible MHD equations.
For stationary full MHD equations, Schneebeli and Sch{\"o}tzau \cite{sch03} proposed a new mixed finite element method
where $\BB$ is discretized by N{\'e}d{\'e}lec's edge elements.
In 2004, Sch{\"o}tzau proved optimal error estimates of the finite element method \cite{sch04}.
In 2010, Greif et al extended the work in \cite{sch04} by discretizing $\Bu$ with $\Hdiv$-conforming face elements
\textit{so that $\Div\Bu_h=0$ holds exactly} \cite{gre10}.
Our choice in this paper for velocity discretization \textit{is exactly the same as} \cite{gre10}.
In 2008, Prohl proved the convergence and error estimates of finite element method for time-dependent MHD equations where
$\BB$ is discretized with $\BH(\curl, \Omega)$-conforming edge elements \cite{pro08}.
In 2017, Hu et al proposed a stable finite element method which discretizes $\BB$ with face elements
and $\BE$ with edge element so that $\Div \BB_h=0$ holds exactly  \cite{hu17}.
In 2018, Hiptmair et al proposed
a new stable finite element method using $(\Bu,p,\BA)$ so that the discrete velocity and the discrete magnetic induction are both divergence-free exactly.
For finite element error estimates, we would also like to mention \cite{he2015IMA} for Euler semi-implicit
scheme and \cite{su2018JSC} for penalty-based finite element methods. We remark that,
for $\BB$-based formulations like in \cite{gun91,sch03, sch04,pro08,he2015IMA} using
edge element or nodal element, the discrete
current density $\BJ_h := \curl \BB_h$ in the moment equation is divergence-free exactly.
We also refer to \cite{ni2012III, zhang2014,li2019SIAM} for charge-conservative methods for inductionless MHD equations
where the applied magnetic induction is known in advance.

The second objective of this paper is to propose a monolithic
iterative solver with augmented block preconditioner
for the discrete problem linearized by Picard's method.
We refer to the monograph of Elman, Silvester, and Wathen
\cite{elman2014Book} for a comprehensive introduction of preconditioners
and iterative solvers for solving incompressible Navier-Stokes equations. For MHD equations, there are extensive studies in the literature, such as \cite{phi16, chacon2008, li17, phillips2014}, on block preconditioners based on approximate Schur complements. We refer to
\cite{shadid2010, shadid2016, lin2018} for algebraic multigrid methods and to \cite{adler2016} for geometric multigrid method.
In \cite{badia2014}, Badia, Mart\'{i}n, and Planas proposed block recursive LU preconditioners for solving thermally coupled inductionless MHD equations.
In the present paper, we follow similar ideas as in \cite{phi16, li17} to propose an augmented Lagrangian block preconditioner for solving the discrete problem.
Since $\Div\Bu_h=0$ in our case, the augmented term $\alpha(\Div\Bu_h,\Div\Bv_h)$ added to the momentum equation does not modify the discrete problem, but enhances the robustness of the block preconditioner.
Numerical examles show that the convergence of preconditioned GMRES solver is quasi-uniform to the number of degrees of freedom (DOFs).

The paper is organized as follows:
In Section 2, we derive a mixed formulation for the full MHD model by means of constrained transport method.
In Section 3, we propose a stable mixed finite element method for the incompressible MHD model
such that $\Bu_h$, $\BJ_h$, and $\BB_h$ are divergence-free exactly in the momentum equations.
A Picard-type linearization for the nonlinear discrete problem is also proposed and the well-posedness of the linearized problem is proven.
In Section 4, to solve the linearized discrete problem, we propose an augmented Lagrangian block preconditioner by deriving the approximate Schur
complements of the MHD system.
In Section 5, we present several numerical examples to verify the theoretical results and to demonstrate the performance of
the preconditioned monolithic solver.
Finally in Section 6, a conclusion is given and some future researches are pointed out.
Throughout the paper, we assume that $\rho, \nu,\mu,\sigma$
are positive constants and denote vector-valued quantities by
boldface notations, such as $\Ltwov:=(\Ltwo)^3$.

\section{Constrained transport model of MHD equations}

Let $\Ltwo$ be the usual Hilbert space of square integrable functions and $\Hone$, $\Hcurl$, $\Hdiv$ be its subspaces with square integrable gradients, curls, and divergences, respectively. Let $\zbHone$, $\zbHcurl$, $\zbHdiv$ denote their subspaces with vanishing traces, vanishing tangential traces, and vanishing normal traces on $\Gamma$ respectively. We refer to \cite[page 26]{gir12} for their definitions and inner products. For convenience, we also introduce the curl-free and divergence-free subspaces
 \begin{align*}
 &\kHcurl=\{\Bv\in\Hcurl:\curl\Bv=\mathbf{0}\},\\
 &\kHdiv=\{\Bv\in\Hdiv:\Div\Bv={0}\}.
 \end{align*}

\subsection{Constrained transport formula}

As remarked in \cite[Section III]{evans1988}, when $\Div\BB_h\ne 0$, the discrete Lorentz force $\BJ_h\times\BB_h$ may yield inaccurate results in the momentum equation of Navier-Stokes equations.
To improve the trustworthiness of $\BB_h$,
one should keep the constraint $\nabla\cdot \BB_h = 0$ in the discrete level.

The traditional CT method is used to solve ideal magnetic induction equations
\begin{equation}\label{CTPoint}
\partial_t \BB + \curl\BE = 0,\qquad
\BE=\BB\times\Bu,\qquad \Div\BB=0.
\end{equation}
In \cite{evans1988}, $\BB$ and $\BE$ are discretized on staggered grids, namely, the discretization of $\BB$ is face-centered and the discretization of $\BE$ is edge-centered.
Since $\Div\BB=0$, there is a magnetic vector potential $\BA$ such that $\BB=\curl\BA$. The CT method amounts to computing the discrete magnetic potential $\BA_h$ by edge-centered discretization and defining the discrete magnetic induction by $\BB_h=\curl\BA_h$.
The evolution equation of $\BA$ is derived from \eqref{CTPoint} by using temporal gauge $\BE=-\partial_t\BA$ (see \cite[(4.14)-(4.17)]{evans1988})
\begin{equation}\label{eq:temp}
\partial_t \BA = \Bu \times \curl\BA.
\end{equation}
Similar ideas can also be found in works related to vector potential methods \cite{hip18,ramshaw1983,rossmanith2006}.

For the stationary MHD model \eqref{sys:MHD}, one faces a major difficulty when dealing with Lorentz force
with only magnetic vector potential (see \cite{evans1988})
 \ben
 \BJ\times\BB = \curl\BH\times \BB
 =\curl(\mu^{-1}\curl\BA) \times \curl\BA,
 \een
namely, one has to discretize second derivatives of $\BA$ in the Lorentz force term.
To overcome the difficulty, we propose to compute the Lorentz force by
 \ben
 \BJ\times\BB = \curl\BH\times \curl\BA
 \een
and discretize $\BH$, $\BA$ as individual variables. This amounts to applying vector potential methods to both $\BJ$ and $\BB$ simultaneously.

\subsection{Constrained transport formula of \eqref{sys:MHD}}

From the vector potential theorem in \cite{amrouche1998}, we easily get the following theorem.

\begin{lemma}\label{lm:VPHcurl}
Suppose $\Omega$ is a simply-connected and Lipschitz domain. For any $\BB\in\kHdiv\cap\zbHcurl$, there exists a unique $\BA\in\Hcurl$ such that
 \begin{subequations}\label{eqs:VP}
 \begin{align}
 \curl\BA=\BB,\quad
 \Div{\BA}=0\quad&\text{in}\;\;\Omega,\\
 \curl\BA\times\Bn=0,\quad
 \BA\cdot\Bn=0\quad&\text{on}\;\;\Gamma.
 \end{align}
 \end{subequations}
\end{lemma}

Because $\BB = \mu\BH$ and taking curls of the first equation of \eqref{eqs:VP} one obtains
\begin{equation}\label{eqs:mainV}
\curl\mu^{-1}\curl\BA = \curl\BH,\quad
\Div{\BA}=0\quad \text{in}\;\;\Omega
\end{equation}
which is precisely the classical double curl problem for $\BA$ and
can be efficiently solved using existing techniques. Using $\BJ = \curl\BH$ and
eliminating the electric field $\BE$ we firstly have
\[\curl(\sigma^{-1}\curl\BH) + \curl(\BB\times\Bu) = \mathbf{0}\quad \text{in}\;\;\Omega\]
where $\BB\times\Bu$ is called induced electric field due to movement of the fluid.
Based on \eqref{eqs:mainV}, then using $\BJ = \curl\BH$ and $\BB = \curl\BA$ both
in the Lorentz force $\BJ\times\BB$ and induced electric field $\BB\times\Bu$, an
equivalent CT form of \eqref{sys:MHD} can be given as follows
\begin{subequations}\label{CTPDEmain}
\begin{align}
\rho\Bu\cdot\nabla\Bu+\nabla{p}-\nu\Delta\Bu
-\curl\BH\times\curl\BA=\Bf\quad&\text{in}\;\;\Omega,\\
\curl(\sigma^{-1}\curl\BH)+\curl(\curl\BA\times\Bu)
=\mathbf{0}\quad&\text{in}\;\;\Omega,\\
\curl\mu^{-1}\curl\BA-\curl\BH=\mathbf{0}\quad
&\text{in}\;\;\Omega,\\
\Div\Bu=0,\quad\Div(\mu\BH)=0,\quad\Div\BA=0
\quad&\text{in}\;\;\Omega,\\
\Bu=\mathbf{0},\quad\BH\times\Bn=\mathbf{0},
\quad\BA\cdot\Bn=0,\quad\curl\BA\times\Bn=\mathbf{0}
\quad&\text{on}\;\;\Gamma.
\end{align}
\end{subequations}
The reason we use $\BA$ to represent the induced electric field is that we want the whole FEM to be stable.
This skill has already appeared in the previous work of \cite{gun91,sch03,sch04}.
Remember that $\rho$, $\nu$, $\sigma$, $\mu$ are all positive constants.
Let $U,H,L$ be the characteristic quantities for velocity, magnetic field, and length of the system respectively.
Define the dimensionless Reynolds number $R_e$, coupling number $\kappa$, and magnetic Reynolds number $R_m$ by
 \ben
 R_e = U L \rho/\nu, \qquad
 \kappa = \mu H^2/(\rho U^2), \qquad
 R_m = \mu\sigma U L.
 \een
Then the MHD system can be written into a dimensionless form
\begin{framed}
\begin{flushleft}
\begin{subequations}\label{eqs:Hmv}
\begin{align}
\Bu\cdot\nabla\Bu+\nabla{p}-R_e^{-1}\Delta\Bu
-\kappa\curl\BH\times\curl\BA=\Bf\quad&\text{in}
\quad\Omega,\label{eq:Hmv1}\\
\kappa R_m^{-1}\curl\curl\BH + \kappa\curl(\curl\BA\times\Bu)+\nabla{r}
=\mathbf{0}\quad&\text{in}\;\;\Omega,\label{eq:Hmv2}\\
\curl\curl\BA-\curl\BH+\nabla\phi=
\mathbf{0}\quad&\text{in}\;\;\Omega,\label{eq:Hmv3}\\
\Div\Bu=0,\quad\Div\BH=0,\quad\Div\BA=0
\quad&\text{in}\;\;\Omega,\label{eq:Hmv4}\\
\BH\times\Bn=\mathbf{0},
\quad\BA\cdot\Bn=0,\quad\curl\BA\times\Bn=\mathbf{0}
\quad&\text{on}\;\;\Gamma, \label{eq:Hmv5}\\
\Bu=\mathbf{0},\quad r=0, \quad \nabla\phi\cdot\Bn =0
\quad&\text{on}\;\;\Gamma. \label{eq:Hmv6}
\end{align}
\end{subequations}
\end{flushleft}
\vspace{-5mm}
\end{framed}

In \eqref{eqs:Hmv}, $r$ and $\phi$ are, respectively, Lagrange multipliers for $\BH$ and $\BA$. Taking divergences
of \eqref{eq:Hmv2} and \eqref{eq:Hmv3} and using \eqref{eq:Hmv6}, we get
\ben
\Delta r= 0, \quad\Delta\phi =0 \quad\mathrm{in}\;\;\Omega,\qquad
r=0, \quad \nabla\phi\cdot\Bn =0  \quad \text{on}\;\;\Gamma.
\een
This implies $\nabla r=\nabla\phi =0$. So \eqref{eqs:Hmv} is actually equivalent to \eqref{CTPDEmain}.
Although the overall structure of the new formulation is similar to the ones in \cite{gun91,sch04,gre10},
here we introduce an extra magnetic vector potential
$\BA$ inspired by original CT methods and therefore ensure the divergence-free
conditions for $\BB_h$ in the Lorentz force.

In our previous work we have used vector potential with edge element method in 3D to represent $\BB$ \cite{hip18}
but the discrete current density of the Lorentz force there is only weakly divergence-free.
In \cite{gre10}, $\Bu_h$ and $\BJ_h$ are divergence-free in the momentum equation.
In the present work we incorporate the two advantages of the formulation in \cite{gun91,gre10} and
traditional vector potential methods to realize the
two divergence-free conditions in the Lorentz force at the same time.
However compared with the work in \cite{gun91,gre10},
two extra variables $\BA$ and $\phi$ are incurred.
One can see that standard mixed finite element methods will give a triple saddle-point problem which is more difficult to solve.
Thus another object of this paper is to develop preconditioned iterative methods to reduce the overall cost as best we can.

\subsection{A weak formulation}
For convenience, we introduce some notations for function spaces
\begin{align*}
\begin{array}{ccc}
\BV:=\zbHonev,\qquad & \BW:=\zbHcurl,\qquad & \BD:=\Hcurl,
 \vspace{1mm}\\
 {Q}:=\Ltwo/\bbR,\qquad & {S}:=\zbHone,\quad
 \qquad &{Y}:=\Hone/\bbR.
\end{array}
\end{align*}
The divergence-free subspaces of $\BV$, $\BW$, and $\BD$ are defined by
 \ben
 \BU(\Div0):=\BU\cap\kHdiv\quad \hbox{for}\;\;
 \BU=\BV,\BW,\BD.
 \een

Multiplying both sides of (\ref{eq:Hmv1}) with $\Bv\in\BV$ and integrating by parts, we get
 \begin{align}\label{eq:IBP1}
 \mathscr{A}(\Bu,\Bv)+\mathscr{O}(\Bu;\Bu,\Bv)
 -\mathscr{L}(\BA;\Bv,\BH) -\SP{p}{\Div\Bv}
 =\SP{\Bf}{\Bv},
 \end{align}
where the bilinear form $\mathscr{A}$ and the trilinear forms $\mathscr{O}$, $\mathscr{L}$ are defined respectively by
 \ben
 &&\mathscr{A}(\Bw,\Bv):=R_e^{-1}\SP{\nabla\Bw}{\nabla\Bv},
 \qquad\mathscr{O}(\Bw;\Bu,\Bv):=\SP{\Bw\cdot\nabla\Bu}{\Bv}, \\
 &&\mathscr{L}(\BA;\Bv,\BH):=\kappa\SP{\curl\BA\times\Bv}{\curl\BH}.
 \een
Multiply both sides of (\ref{eq:Hmv2}) with $\Bw\in\BW$ and both sides of (\ref{eq:Hmv3}) with $\Bd\in\BD$.
Using integration by parts, we get
 \begin{align}
 &\mathscr{C}_1(\BH,\Bw)+
 \mathscr{L}(\BA;\Bu,\Bw)+\SP{\nabla{r}}{\Bw}=0, \label{eq:IBP2}\\
 &\mathscr{C}(\BA,\Bd)-\SP{\BH}{\curl\Bd}+\SP{\nabla\phi}{\Bd}=0,
 \label{eq:IBP3}
 \end{align}
where the bilinear forms are defined by
 \ben
 \mathscr{C}_1(\Bv,\Bw):=\kappa{R}_m^{-1}(\curl\Bv,\curl\Bw),
 \qquad
 \mathscr{C}(\Bv,\Bw):=(\curl\Bv,\curl\Bw).
 \een
Combining (\ref{eq:IBP1})--(\ref{eq:IBP3}), we obtain a weak formulation of (\ref{eqs:Hmv}):
\begin{framed}
\begin{flushleft}
Find $(\Bu,\BH,\BA)\in\BV\times\BW\times\BD$ and $(p,r,\phi)\in{Q}\times{S}\times{Y}$ such that
 \begin{subequations}\label{eqs:weak}
 \begin{align}
 &\mathscr{A}(\Bu,\Bv)+\mathscr{O}(\Bu;\Bu,\Bv)
 -\mathscr{L}(\BA;\Bv,\BH)-\SP{p}{\Div\Bv}
 =\SP{\Bf}{\Bv},\label{eq:weak1}\\
 &\mathscr{C}_1(\BH,\Bw)
 +\mathscr{L}(\BA;\Bu,\Bw)
 +\SP{\nabla{r}}{\Bw}=0,\label{eq:weak2}\\
 &\mathscr{C}(\BA,\Bd)-\SP{\BH}{\curl\Bd}
 +\SP{\nabla\phi}{\Bd}=0,
 \label{eq:weak3}\\
 &\SP{\Div\Bu}{q}=0,\qquad\SP{\BH}{\nabla{s}}=0,
 \qquad\SP{\BA}{\nabla\varphi}=0,\label{eq:weak4}
 \end{align}
 \end{subequations}
for all $(\Bv,\Bw,\Bd)\in\BV\times\BW\times\BD$ and $(q,s,\varphi)\in Q\times S\times Y$.
\end{flushleft}
\end{framed}

\begin{theorem}\label{thm:reduceB}
The solutions of \eqref{eqs:weak} satisfy the stability estimate
\begin{align}\label{ineq:reducedB}
\NHonev{\Bu}+\NHcurl{\BH}+\NHcurl{\BA}\le C\NLtwov{\Bf},
\end{align}
where the constant $C$ depends only on $\kappa$, $R_e$, $R_m$, and the domain $\Omega$.
\end{theorem}
\begin{proof}
Using \eqref{eq:weak4}, it is easy to see $\Bu\in\BV(\Div0)$, $\BH\in\BW(\Div0)$, and $\BA\in\BD(\Div0)\cap\zbHdiv$. Then
\eqref{eqs:weak} is reduced to
\begin{align*}
 &\mathscr{A}(\Bu,\Bv) +\mathscr{O}(\Bu;\Bu,\Bv)
 -\mathscr{L}(\BA;\Bv,\BH)
 =\SP{\Bf}{\Bv} \qquad \forall\,\Bv\in \BV(\Div0),\\
 &\mathscr{C}_1(\BH,\Bw)+\mathscr{L}(\BA;\Bu,\Bw)=0
  \qquad \forall\,\Bw\in \BW(\Div0),\\
 &\mathscr{C}(\BA,\Bd)-\SP{\BH}{\curl\Bd}=0
  \qquad \forall\,\Bd\in \BD(\Div0).
 \end{align*}
Taking $(\Bv,\Bw,\Bd)=(\Bu,\BH,\BA)$ and using $\Div\Bu=0$, we find that
\begin{align*}
&\frac{1}{R_e}\SNHonev{\Bu}^2-\mathscr{L}(\BA;\Bu,\BH)
  =\SP{\Bf}{\Bu},\qquad
\frac{\kappa}{R_m}\NLtwov{\curl\BH}^2
 +\mathscr{L}(\BA;\Bu,\BH)=0,\\
&\NLtwov{\curl\BA}^2 =\SP{\BH}{\curl\BA}.
 \end{align*}
Adding up the first and second equations yields
\begin{align*}
\frac{1}{R_e}\SNHonev{\Bu}^2
+\frac{\kappa}{R_m}\NLtwov{\curl\BH}^2=(\Bf,\Bu)
\le C\NLtwov{\Bf}^2 + \frac{1}{2R_e}\SNHonev{\Bu}^2,
\end{align*}
that is, $\NHonev{\Bu}
+\NLtwov{\curl\BH}\le C\NLtwov{\Bf}$, where we have used Poincar\'{e}'s inequality to $\Bu$.
Moreover, applying Poincar\'{e}-type inequality to $\BH\in\BW(\Div0)$ and $\BA\in\BD(\Div0)\cap\zbHdiv$ (cf. e.g. \cite{amrouche1998}), we obtain
 \ben
 &&\NHcurl{\BH}\le C\left(\NLtwov{\curl\BH}+\NLtwo{\Div\BH}\right)
 =C\NLtwov{\curl\BH},\\
 &&\NHcurl{\BA}\le C\NLtwov{\curl\BA}
 \le C\NLtwov{\BH}\le C\NLtwov{\curl\BH}.
 \een
The proof is completed.
\end{proof}

\section{Mixed finite element approximation}
In this section, we study finite element approximation to the weak formulation of the MHD model.
Inspired by \cite{gre10,hip18}, the velocity $\Bu$ will be discretizd by $\Hdiv$-conforming
Brezzi-Douglas-Marini (BDM) elements and a DG-type formulation with interior penalties.
Let $\mesh$ be a quasi-uniform and shape-regular tetrahedral mesh of $\Omega$.
Let $h_K$ be the diameter of a tetrahedron $K\in\Ct_h$ and
let $h=\max\limits_{K\in\mesh}h_K$ denote the mesh size of $\mesh$.

\subsection{An interior-penalty finite element method}
First we introduce the finite element spaces for $\Bu$, $\BH$, and $\BA$ as follows
 \begin{align*}
 \BV_h&:=\{\Bv\in\zbHdiv:\;\;\Bv|_K\in\BP_1(K),\;\;\;\forall{K}\in\mesh\},\\
 \BW_h&:=\{\Bw\in\BW:\;\;\Bw|_K\in\BP_1(K),\;\;\;\forall{K}\in\mesh\},\\
 \BD_h&:=\{\Bd\in\BD:\;\;\;\Bd|_K\in\BP_1(K),\;\;\;
 \forall{K}\in\mesh\},
 \end{align*}
where $\BP_k(K)=(P_k(K))^3$ and $P_k$ is the space of polynomials with degree $\le{k}$. Functions in $\BV_h$ are continuous normally but may be discontinuous  tangentially, while functions in $\BW_h\cup\BD_h$ are continuous tangentially but may be discontinuous normally. The finite element spaces for multipliers $(p,r,\phi)$ are defined respectively by
 \begin{align*}
 Q_h&:=\{q\in{Q}:\;\;q|_K\in{P}_0(K),\;\;\;\forall{K}\in\mesh\},\\
 S_h&:=\{v\in{S}:\;\;v|_K\in{P}_2(K),\;\;\;\forall{K}\in\mesh\},\\
 Y_h&:=\{\varphi\in Y:\;\;\varphi|_K\in{P}_2(K),\;\;\;
 \forall{K}\in\mesh\}.
 \end{align*}
Clearly functions in $S_h$ and $Y_h$ are continuous.

Let $\Cf_h$ denote the set of faces of all tetrahedra in $\Ct_h$.
We endow each $F\in\Cf_h$ with a unit normal $\Bn_F$ which points to the exterior of $\Omega$ when $F\subset\Gamma$ and to $K_-$ when
$F=\partial{K}_+\cap\partial{K}_-$ for two adjacent elements $K_\pm\in\Ct_h$.
Let $\varphi$ be a scalar-, vectorial, or matrix-valued function which is piecewise smooth over $\mesh$.
The mean value and the jump of $\varphi$ on $F$ are defined respectively by
 \ben
 \{\{\varphi\}\}:=(\varphi_++\varphi_-)/2,\quad
 \jump{\varphi}:=\varphi_+-\varphi_-\quad\text{on}\;\;{F},
 \een
where $\varphi_{\pm}$ denote the traces of $\varphi$ on $F$ from inside of $K_{\pm}$ respectively.
For any face $F=\partial{K}_+\cap\Gamma$, the mean value and the jump of $\varphi$ on $F$ are defined by
 \ben
 \{\{\varphi\}\}=\jump{\varphi}=\varphi_+\quad\text{on}\;\;{F}.
 \een

The discrete counterparts of $\mathscr{A}$ and $\mathscr{O}$ are defined by
\begin{align}
 \mathscr{A}_h(\Bu,\Bv)=
 &\frac{1}{R_e}\sum_{K\in\mesh}\int_K \nabla\Bu:\nabla\Bv+\frac{\gamma}{R_e} \sum_{F\in\Cf_h}h_F^{-1}\int_F\jump{\Bu}\cdot\jump{\Bv}
 \notag\\
 &-\frac{1}{R_e}\sum_{F\in\Cf_h}\int_F \left(\{\{\frac{\partial\Bu}{\partial\Bn_F}\}\}
 \cdot\jump{\Bv}+\{\{\frac{\partial\Bv}{\partial\Bn_F}\}\}\cdot\jump{\Bu}\right),\label{eq:Ah}\\
 \mathscr{O}_h(\Bw;\Bu,\Bv)=
 &-\sum_{K\in\mesh}\int_K\Bu\cdot
 \Div(\Bw\otimes\Bv)+\sum_{K\in\mesh}
 \int_{\partial{K}}(\Bw\cdot\Bn_K)(\Bu^\downarrow\cdot\Bv),
 \label{eq:Oh}
\end{align}
where $\gamma>0$ is the penalty parameter, $\Bn_K$ the unit outer normal of $\partial{K}$,
$h_F$ the diameter of $F$, and $\Bu^\downarrow$ the upwind convective flux defined by
 \ben
 \Bu^\downarrow(\Bx)=
 \begin{cases}
 \lim\limits_{\epsilon\to0^+}\Bu(\Bx-\epsilon\Bw(\Bx)),
 \quad&\Bx\in\partial{K}/\Gamma,\\
 \mathbf{0},\quad&\Bx\in\partial{K}\cap\Gamma.
 \end{cases}
 \een
In \eqref{eq:Ah}, the second term in the right-hand side represents interior
penalties used to insure the stability of discrete solutions.

The finite element approximation to problem \eqref{eqs:weak} is given as follows:
\begin{framed}
\begin{flushleft}
Find $(\Bu_h,\BH_h,\BA_h)\in\BV_h\times\BW_h\times\BD_h$ and $(p_h,r_h,\phi_h)\in{Q}_h\times{S}_h\times{Y}_h$ such that
 \begin{subequations}\label{weakh}
 \begin{align}
 &\mathscr{A}_h(\Bu_h,\Bv_h)+ \mathscr{O}_h(\Bu_h;\Bu_h,\Bv_h)
 -\mathscr{L}(\BA_h;\Bv_h,\BH_h)
 -\SP{p_h}{\Div\Bv_h} =\SP{\Bf}{\Bv_h},\label{eq:weakh1}\\
 &\mathscr{C}_1(\BH_h,\Bw_h) +\mathscr{L}(\BA_h;\Bu_h,\Bw_h)
 +\SP{\nabla{r}_h}{\Bw_h}=0,\label{eq:weakh2}\\
 &\mathscr{C}(\BA_h,\Bd_h)-\SP{\BH_h}{\curl\Bd_h}
 +\SP{\nabla\phi_h}{\Bd_h}=0,
 \label{eq:weakh3}\\
 &\SP{\Div\Bu_h}{q_h}=0,\qquad\SP{\BH_h}{\nabla{s}_h}=0,
 \qquad\SP{\BA_h}{\nabla\varphi_h}=0,\label{eq:weakh4}
 \end{align}
 \end{subequations}
for all $(\Bv_h,\Bw_h,\Bd_h)\in\BV_h\times\BW_h\times\BD_h$ and all $(q_h,s_h,\varphi_h)\in Q_h\times S_h\times Y_h$.
\end{flushleft}
\end{framed}
Due to \eqref{eq:weakh4}, we have $\SP{\Div\Bu_h}{q_h}=0$
for any $q_h \in Q_h$. Because $\Div\Bu_h \in Q_h$ also holds due
to the definition of finite element space, thus we obtain
\begin{equation}
(\Div\Bu_h,\Div\Bu_h) = 0 \Rightarrow \|\Div\Bu_h\|_{L^2} = 0
\end{equation}
Because $\Bu_h \in \BH(\Div,\Omega)$, we have $\Div\Bu_h = 0$ exactly on the discrete level.
This assertion directly comes from the work in \cite{gre10}.
However, here we would like to mention that in the work \cite{hu17},
the $\BB$ is also discretized by $\Hdiv$-conforming
element (Raviart-Thommas element there which is similar to BDM element) and
one will only have $\|\Div\BB_h\|_{L^2}$ = 0 (see Lemma 2 of \cite{hu17}).
The authors in \cite{hu17} claim that $\nabla\cdot\BB_h = 0$ holds exactly on the discrete level, so
one can similar assertion that $\Div\Bu_h = 0$ holds exactly. The main difference between our CT-FEM methods
and the $\BB-\BE$ methods in \cite{hu17} is that we use magnetic potential $\BA$ and magnetic field $\BH$
as independent variables and therefore the divergence-free conditions for $\BJ$ and $\BB$ in the Lorentz force
are both satisfied naturally. Moreover as indicated by Theorem \ref{thm:picardB} in the following, our new finite
element method is also energy stable.

We also remark that if $\mu$ is constant,
the scheme \eqref{weakh} only enforces the constraint of $\BH_h$ weakly.
But our main concern is the divergence-free constraints in the momentum equations,
especially in the Lorentz force $\BJ_h\times\BB_h$,
so we compromise and relax $\BH_h$'s constraint and
just adopt edge element to discretize $\BH_h$.

To make the $\BB_h$ in the Lorentz force divergence-free,
we recover a discrete magnetic vector potential $\BA_h$
by solving a double curl problem using edge finite element method,
\begin{equation}
\mathscr{C}(\BA_h,\Bd_h)+\SP{\nabla\phi_h}{\Bd_h}=\SP{\BH_h}{\curl\Bd_h},
\quad \SP{\BA_h}{\nabla\varphi_h}=0
\end{equation}
and then use $\BA_h$ to compute discrete magnetic induction $\BB_h=\curl\BA_h$ in the Lorentz force.
This philosophy has been used successfully in projection methods \cite{brackbill1980} and
high-order CT methods on unstaggered grid \cite{rossmanith2006}.
In \cite[Section 4.2]{rossmanith2006}, for time-dependent ideal compressible MHD equaitons,
Rossmanith solved the MHD equations using a "base scheme" to firstly obtain an intermediate magnetic field,
which is not divergence-free. And then Rossmanith use the intermediate magnetic field to evolve the magnetic
potentials to obtain the precisely divergence-free magnetic filed in the next time-step.
We remark that $\BH_h$ in our methods plays the same role as
the pre-computed magnetic filed in Rossmanith's methods and other CT methods.
The difference is that we solve a double curl problem instead to reconstruct the $\BB_h$ in the Lorentz force.
And this is the reason we call our methods by constrained transport divergence-free finite element.

\subsection{An iterative scheme for the discrete problem}
The discrete problem \eqref{weakh} is a nonlinear system. Here we propose an iterative scheme of Picard type to solve the problem.
Let $\Bu_h^{n-1}\in\BV_h$, $\BA_h^{n-1}\in\BD_h$, $n\ge 1$, be the approximate solutions in the $(n-1)^{\rm th}$ iteration.
The approximate solutions in the $n^{\rm th}$ iteration solve the coupled linear system:
\begin{framed}
\begin{flushleft}
Find $(\Bu_h^n,\BH_h^n,\BA_h^n)\in\BV_h\times\BW_h\times\BD_h$ and $(p_h^n,r_h^n,\phi_h^n)\in{Q}_h\times{S}_h\times{Y}_h$ such that
 \begin{subequations}\label{weakL}
 \begin{align}
 &\mathscr{A}_h(\Bu_h^n,\Bv_h)+
 \mathscr{O}_h(\Bu_h^{n-1};\Bu_h^n,\Bv_h)
 -\mathscr{L}(\BA_h^{n-1};\Bv_h,\BH_h^n)
 -\SP{p_h^n}{\Div\Bv_h} =\SP{\Bf}{\Bv_h},\label{eq:weakL1}\\
 &\mathscr{C}_1(\BH_h^n,\Bw_h)
 +\mathscr{L}(\BA_h^{n-1};\Bu_h^n,\Bw_h)
 +\SP{\nabla{r}^n_h}{\Bw_h}=0,\label{eq:weakL2}\\
 &\mathscr{C}(\BA_h^n,\Bd_h)-\SP{\BH_h^n}{\curl\Bd_h}
 +\SP{\nabla\phi_h^n}{\Bd_h}=0,
 \label{eq:weakL3}\\
 &\SP{\Div\Bu_h^n}{q_h}=0,\qquad
 \SP{\BH_h}{\nabla{s}^n_h}=0,
 \qquad\SP{\BA_h^n}{\nabla\varphi_h}=0,\label{eq:weakL4}
 \end{align}
 \end{subequations}
for all $(\Bv_h,\Bw_h,\Bd_h)\in\BV_h\times\BW_h\times\BD_h$ and all $(q_h,s_h,\varphi_h)\in Q_h\times S_h\times Y_h$.
\end{flushleft}
\end{framed}
The implicit upwind DG term $\mathscr{O}_h(\Bw_h;\Bu_h,\Bv_h)$ is
tailored to the standard convection term $\Bw\cdot\nabla\Bu$,
which is difficult to modify for Newton's method. Also see Remark 3.3 of \cite{gre10}. When upwinding is not used,
Newton¡¯s method can be straightforwardly applied, which is just the case $\BP_2-P_1$ Taylor-Hood element for $\Bu_h-p_h$
pair \cite{li17}. We admit that this is really a drawback for stationary problems where Newton's iteration may be more efficent
nonlinear solver with relative large parameters (see numerical examples in \cite{li17}). The numerical experiments in the following
also indicate that Picard iteration is not robust enough for large $R_e$ and $R_m$.
In the future, \textit{acceleration techniques in optimization methods} for nonlinear iteration
can be incorporated for the present discretization.

\subsection{Well-posedness of \eqref{weakL}}

First we introduce the discrete semi-norm and norm for piecewise regular functions
 \begin{equation*}
 \Nsemi{\varphi}:=\left(\sum_{K\in\mesh}
 \SNHone[K]{\varphi}^2\right)^{1/2},\quad
 \Ndg{\varphi}:=\left(\Nsemi{\varphi}^2
 +\sum_{F\in\Cf_h}h_F^{-1}\int_F\jump{\varphi}^2\right)^{1/2}.
 \end{equation*}
The following lemma states that $\mathscr{A}_h$ is coercive and continuous and $\mathscr{O}_h$ is positive and continuous.

\begin{lemma}{\rm\cite[Lemma 3.2--3.3]{hip18}}\label{lm:AhCh}
Suppose $\gamma$ is large enough but independent of $h_F$ and $R_e$. There are constants $\theta_1,\theta_2>0$ independent of $h_F$ and $R_e$ such that
 \ben
 \mathscr{A}_h(\Bv,\Bv)\ge \theta_1R_e^{-1}\Ndg{\Bv}^2,
 \quad\mathscr{A}_h(\Bu,\Bv)\le \theta_2R_e^{-1}
 \Ndg{\Bu}\Ndg{\Bv}\quad\forall\,\Bu,\Bv\in\BD_k(\mesh),
 \een
where $\BD_k(\mesh)=\{\Bv\in\Ltwov: \Bv|_K\in\BP_k(K),\;\forall\, K\in\Ct_h\}$. Moreover, there is a constant $C$ independent of $h$ such that, for any $\Bw,\Bw_1\in\BD_k(\mesh)\cap\kHdiv$,
 \ben
 &&\mathscr{O}_h(\Bw;\Bv,\Bv)=\frac12\sum_{F\in\Cf_h}
 \int_F\SN{\Bw\cdot\Bn}\SN{\jump{\Bv}}^2,\\
 &&\SN{\mathscr{O}_h(\Bw;\Bu,\Bv)-\mathscr{O}_h(\Bw_1;\Bu,\Bv)}
 \le{C}\Ndg{\Bw-\Bw_1}\Ndg{\Bu}\Ndg{\Bv}.
 \een
\end{lemma}

For convenience, we define the divergence-free subspace of $\BV_h$ by
 \ben
 \BV_h(\Div0)=\BV_h\cap\kHdiv,
 \een
and define the weakly divergence-free subspaces of $\BW_h,\BD_h$ by
 \begin{align*}
 \BW_h(\Div0)&:=\{\Bc_h\in\BW_h:\;\;{\SP{\Bc_h}{\nabla{s}_h}=0},
 \;\;\;\forall{s}_h\in{S}_h\},\\
 \BD_h(\Div0)&:=\{\Bd_h\in\BD_h:\;\;
 {\SP{\Bd_h}{\nabla\varphi_h}=0},\;\;\;
 \forall\varphi_h\in{Y}_h\}.
 \end{align*}
Generally, we have $\BW_h(\Div0)\not\subset\BW(\Div0)$ and $\BD_h(\Div0)\not\subset\BD(\Div0)$. However, the $\Ltwov$-orthogonal decompositions or discrete Helmholtz decompositions hold
 \begin{align}\label{Ddecomp}
 \BW_h=\BW_h(\Div0)\oplus\nabla{S}_h,\qquad
 \BD_h=\BD_h(\Div0)\oplus\nabla{Y}_h.
 \end{align}
From \cite[Theorem 4.7]{hip02}, there is a constant $C$ independent of $h$ such that the discrete Poincar{\'{e}} inequality holds
 \begin{align}\label{ineq:disL2}
 \NLtwov{\Bc_h}\le{C}\NLtwov{\curl\Bc_h}\qquad
 \forall\,\Bc_h\in\BW_h(\Div0)\cup\BD_h(\Div0).
 \end{align}

From \eqref{eq:weakL4}, we easily find that
 \ben
 \Bu_h^n\in \BV_h(\Div0),\qquad
 \BH_h^n\in\BW_h(\Div0),\qquad
 \BA_h^n\in\BD_h(\Div0).
 \een
So \eqref{weakL} can be written into a reduced form:
Find $(\Bu_h^n,\BH_h^n)\in\BV_h(\Div0) \times\BW_h(\Div0)$ and $\BA_h^n\in\BD_h(\Div0)$ such that
 \begin{subequations}\label{weakr}
 \begin{align}
 &a((\Bu_h^n,\BH_h^n),(\Bv_h,\Bw_h)) =\SP{\Bf}{\Bv_h}\quad
 \forall\,(\Bv_h,\Bw_h)\in\BV_h(\Div0) \times\BW_h(\Div0), \label{eq:weakr1}\\
 &\mathscr{C}(\BA_h^n,\Bd_h)=\SP{\BH_h^n}{\curl\Bd_h}
 \quad \forall\,\Bd_h\in\BD_h(\Div0),
 \label{eq:weakr2}
 \end{align}
 \end{subequations}
where the bilinear form $a$ is defined by
 \begin{align*}
 a((\Bu_h^n,\BH_h^n),(\Bv_h,\Bw_h)):=\,
 &\mathscr{A}_h(\Bu_h^n,\Bv_h)
 +\mathscr{O}_h(\Bu_h^{n-1};\Bu_h^n,\Bv_h)
 +\mathscr{C}_1(\BH_h^n,\Bw_h) \\
 &-\mathscr{L}(\BA_h^{n-1};\Bv_h,\BH_h^n)
 +\mathscr{L}(\BA_h^{n-1};\Bu_h^n,\Bw_h).
 \end{align*}

\begin{theorem}\label{thm:picardB}
Problem \eqref{weakL} has unique solutions. There exists a constant $C$ independent of $h$ such that
 \begin{align}\label{ineq:picardB}
 \Ndg{\Bu_h^n}+\NHcurl{\BH_h^n}+ \NHcurl{\BA_h^n}\le C\NLtwov{\Bf}.
 \end{align}
\end{theorem}
\begin{proof}
Since $\Bu_n^{n-1}\in\BV_h(\Div0)$, we have $\Div\Bu_h^{n-1}=0$. From Lemma~\ref{lm:AhCh} and inequality \eqref{ineq:disL2}, it is easy to see that
 \ben
 a((\Bv_h,\Bw_h),(\Bv_h,\Bw_h))
 &\ge& \mathscr{A}_h(\Bv_h,\Bv_h)
 +\kappa R_m^{-1}\NLtwov{\curl\Bw_h}^2 \\
 &\ge& C\left(\N{\Bv}^2_{1,h}+\NHcurl{\Bw_h}\right),
 \een
for all $(\Bv_h,\Bw_h)\in\BV_h(\Div0)\times\BW_h(\Div0)$, where  $C>0$ is a constant independent of $h$. Therefore, the bilinear form $a$ is coercive on $\BV_h(\Div0)\times\BW_h(\Div0)$. So the finite dimensional problem \eqref{eq:weakr1} has unique solutions
$(\Bu_h^n,\BH_h^n)$. Moreover, \eqref{ineq:disL2} implies that $\mathscr{C}$ is coercive on $\BD_h(\Div0)$. So for the $\BH_h^n$ obtained from \eqref{eq:weakr1}, problem \eqref{eq:weakr2}
has a unique solution $\BA_h^n$.

From \cite[Proposition 3.3]{cockburn2007}, there is a constant $C_{\rm inf}>0$ independent of $h$ such that
 \ben
 \sup_{0\ne\Bv_h\in\BV_h}\frac{(\Div\Bv_h,q_h)}{\N{\Bv_h}_{1,h}}
 \ge C_{\rm inf}\NLtwo{q_h}\qquad
 \forall\,q_h\in Q_h.
 \een
Moreover, by the inclusions $\nabla S_h\subset\BW_h$, $\nabla Y_h\subset \BD_h$ and Poincar\'{e}'s inequality, there is a constant $\hat{C}_{\rm inf}>0$ depending only on $\Omega$ such that
 \ben
 \sup_{0\ne\Bw_h\in\BW_h}
 \frac{(\Bw_h,\nabla s_h)}{\NHcurl{\Bw_h}}
 &\ge& \NLtwov{\nabla s_h} \ge \hat{C}_{\rm inf}\NHonev{s_h}
 \qquad \forall\,s_h\in S_h, \\
 \sup_{0\ne\Bd_h\in\BD_h}
 \frac{(\Bd_h,\nabla \varphi_h)}{\NHcurl{\Bd_h}}
 &\ge& \NLtwov{\nabla \varphi_h}
 \ge \hat{C}_{\rm inf}\NHonev{\varphi_h}
 \qquad \forall\,\varphi_h\in Y_h.
 \een
We conclude the existence and uniqueness of $p_h^n$, $r_h^n$, and $\phi_h^n$ from \eqref{eq:weakL1}--\eqref{eq:weakL3}.

Finally, the stability in \eqref{ineq:picardB} can be proven by arguments similar to the proof of Theorem~\ref{thm:reduceB}. We do not elaborate on the details.
\end{proof}
\vspace{2mm}

Since $\BV_h$, $\BW_h$, and $\BD_h$ are finite dimensional, the stability \eqref{ineq:picardB} implies that, upon an extracted subsequence, the linearized solutions converge strongly to three functions
$\Bu_h\in\BV_h(\Div0)$, $\BH_h\in\BW_h(\Div0)$, and $\BA_h\in\BD_h(\Div0)$. Moreover, the limits solve the nonlinear problem \eqref{weakh} and satisfy
 \ben
 \Ndg{\Bu_h}+\NLtwov{\curl\BH_h}+\NLtwov{\curl\BA_h}
 \le{C}\NLtwov{\Bf}.
 \een
The convergence of the original sequence $\{(\Bu_h^n,\BH_h^n,\BA_h^n)\}$ and the uniqueness of solutions to \eqref{weakh} can be proven by arguments similar to \cite{sch04} upon assuming that $R_e,R_m,\kappa$ are small enough. Since we are only interested in proposing the conservative scheme and its discrete solver, these are beyond the scope of this paper.
\vspace{2mm}

\begin{remark}\label{rm:free}
Because $\BH_h$ and $\BA_h$ belong to $\BH(\curl,\Omega)$, we
have $\BJ_h:=\curl\BH_h\in\BH(\Div,\Omega)$ and $\BB_h:=\curl\BA_h\in\BH(\Div,\Omega)$ \cite{monk03}.
Thus the discrete current density and magnetic induction in the
Lorentz force are precisely divergence-free because we can
directly use divergence operator on $\BJ_h$ and $\BB_h$.
Moreover due to the mixed finite element for $\Bu_h-p_h$ with $\BH(\Div)$-conforming
element for $\Bu_h$ the scheme is
also mass-conservative (see \cite{cockburn2007,john2017,gre10} for details).
Theorem \ref{thm:picardB} indicates that our finite element method is also energy stable.
\end{remark}

\section{An augmented Lagrangian block preconditioner}
The purpose of this section is to propose a preconditioner for the linear algebraic systems
resulting from the the Picard iteration \eqref{weakL}. Because we use a monolithic way and the linear algebraic systems
are a series of triple saddle-point problems, the usual Krylov subspace methods such as GMRES will be extremely slow to converge
without preconditioning.
And it would be ideal if the number of iterations using fixed tolerance did not grow under mesh refinement.
Here we develop an augmented Lagrangian block preconditioner which follows the work in \cite{li17}, where
a grad-div stabilized formulation for the model in \cite{sch04} is used.
The basic ideas for approximate block factorization
and operators' commutativity come from the work in \cite{phillips2014, phi16}.
Note that in our previous work \cite{li17} $\BP_2-P_1$ Taylor-Hood Element is used
for velocity-pressure pair. But now we adopt $\BH(\Div,\Omega)$-conforming element for velocity.
To devise our preconditioner, we shall follow the approximate Schur complement techniques which has already gained much success in
incompressible Navier-Stokes equations \cite{elman2014Book}.

\subsection{Algebraic form of problem~\eqref{weakL}}

Since $\Div\Bu_h^n=0$, we add a grad-div stabilization term in the momentum equation
and rewrite \eqref{weakL} as follows
 \begin{align*}
 &\mathscr{C}(\BA_h^n,\Bd_h)-\SP{\BH_h^n}{\curl\Bd_h}+\SP{\nabla\phi_h^n}{\Bd_h}=0, \\
 &\SP{\BA_h^n}{\nabla\varphi_h}=0,\\
 &\mathscr{C}_1(\BH_h^n,\Bw_h)+\mathscr{L}(\BA_h^{n-1};\Bu_h^n,\Bw_h)+\SP{\nabla{r}^n_h}{\Bw_h}=0, \\
 &\SP{\BH_h}{\nabla{s}^n_h}=0,\\
 &\mathscr{A}_1(\Bu_h^n,\Bv_h)
 +\mathscr{O}_h(\Bu_h^{n-1};\Bu_h^n,\Bv_h)-\mathscr{L}(\BA_h^{n-1};\Bv_h,\BH_h^n)
 -\SP{p_h^n}{\Div\Bv_h} =\SP{\Bf}{\Bv_h}, \\
 &\SP{\Div\Bu_h^n}{q_h}=0,
\end{align*}
for all $(\Bv_h,\Bw_h,\Bd_h)\in\BV_h\times\BW_h\times\BD_h$ and all $(q_h,s_h,\varphi_h)\in Q_h\times S_h\times Y_h$, where
\ben
\mathscr{A}_1(\Bw_h,\Bv_h):= \mathscr{A}_h(\Bw_h,\Bv_h)+\alpha\SP{\Div\Bw_h}{\Div\Bv_h}.
\een
We have rearranged the variables order as $(\BA_h, \phi_h, \BH_h, r_h, \Bu_h, p_h)$
for easy block factorization and preconditioning.
The adding term $\alpha\SP{\Div\Bu_n}{\Div\Bv}$ is called grad-div stabilization \cite{olshanskii2004}
and $\alpha$ is called the grad-div parameter.
Because $\Div\Bu_h = 0$ can be guaranteed \cite{gre10},
it does not change the discrete solutions compared with that using Taylor-Hood element,
but enhances the performance of the preconditioner.
Here we reiterate that though we realize the three divergence-free constraints for the momentum equation,
it is at the expense of more degrees of freedom because we add variables $\BA$ and $\phi$.
The new CT formulation can also be regarded as a modification of the model in \cite{gre10},
namely an extra double curl problem for $\BA$ is solved to
ensure the divergence-free $\BB_h$ in the Lorentz force.
This fact means that for the proposed scheme existing FEM codes can be reused in a simple way
other than a classical double curl problem solve for magnetic vector potential.

The linear system can be written into an algebraic form
\begin{align}\label{eq:algeb}
\bbA\Vx=\Vb,
\end{align}
where $\bbA$ is the stiffness matrix, $\Vx$ the vector of DOFs, and $\Vb$ the load vector. They are given in block forms by
\begin{align*}
\bbA =
\left(
\begin{array}{cccccc}
\bbC & \bbG^\top & \bbK  & 0         & 0         & 0 	\\
\bbG & 0         & 0     & 0         & 0         & 0 	\\
0    & 0         & \bbH  & \bbD^\top & \bbJ^\top & 0 	\\
0    & 0         & \bbD  & 0         & 0         & 0 	\\
0    & 0         & -\bbJ & 0         & \bbF      & \bbB^\top \\
0    & 0         & 0     & 0         & \bbB      & 0
\end{array}
\right),
\quad\Vx=
\left(
\begin{array}{c}
\Vx_A \\
\Vx_\phi \\
\Vx_H \\
\Vx_r \\
\Vx_u \\
\Vx_p
\end{array}
\right),
\quad\Vb=
\left(
\begin{array}{c}
\Vb_A \\
\Vb_\phi \\
\Vb_H \\
\Vb_r \\
\Vb_u \\
\Vb_p
\end{array}
\right).
\end{align*}
Here $\Vx_A$, $\Vx_\phi$, $\Vx_H$, $\Vx_r$, $\Vx_u$, $\Vx_p$ are vectors of DOFs belonging to $\BA_n$, $\phi_n$, $\BH_n$, $r_n$, $\Bu_n$, $p_n$ respectively and $\Vb_A$, $\Vb_\phi$, $\Vb_H$, $\Vb_r$, $\Vb_u$, $\Vb_p$ are the corresponding load vectors. The block matrices $\bbC$, $\bbG$, $\bbK$, $\bbH$,
$\bbD$, $\bbJ$, $\bbF$ and $\bbB$ are Galerkin matrices defined by
\begin{align*}
\begin{array}{ll}
\bbC\leftrightarrow\mathscr{C}(\BA_h^n,\Bd_h),  \qquad \;
\bbG\leftrightarrow\SP{\BA_h^n}{\nabla\varphi_h},
&\quad\bbK\leftrightarrow-\SP{\BH_h^n}{\curl\Bd_h},
    \vspace{1mm}\\
\bbH\leftrightarrow \mathscr{C}_1(\BH_h^n,\Bd_h), \qquad
\bbD\leftrightarrow\SP{\BH_h^n}{\nabla s_h},
&\quad\bbJ\leftrightarrow
\mathscr{L}(\BA_h^{n-1};\Bv_h,\BH_h^n),
    \vspace{1mm}\\
\bbF\leftrightarrow \mathscr{A}_1(\Bu_h^n,\Bv_h)
+\mathscr{O}_h(\Bu_h^{n-1};\Bu_h^n,\Bv_h),
&\quad\bbB\leftrightarrow-\SP{p_h^n}{\Div\Bv_h}.
\end{array}
\end{align*}
We remark that compared with the formulation in \cite{gre10,phi16,li17},
the extra cost for the linear algebraic equations is only the classical double curl saddle problem
for the magnetic vector potential $\BA_h$
\[\left(
    \begin{array}{cc}
      \bbC & \bbG^\top \\
      \bbG & 0 \\
    \end{array}
  \right)
\]

\subsection{Ideal preconditioner of $\bbA$}

Now we deduce a preconditioner of $\bbA$ based on LU factorization and approximate Schur complements.
Since $\bbC$ and $\bbH$ are discretization of $\curl\curl$ operators, they are singular matrices.
We first apply mass augmentation techniques to the two saddle structures for $\bbC$ and $\bbH$ (see more details in \cite{greif2007} for the augmentation of Maxwell saddle-point problem and Section 3.1 of \cite{phi16} for MHD)
\begin{equation}
\left(
  \begin{array}{cc}
    \bbC & \bbG^\top \\
    \bbG & 0 \\
  \end{array}
\right),\quad
\left(
  \begin{array}{cc}
    \bbH & \bbD^\top \\
    \bbD & 0 \\
  \end{array}
\right).
\end{equation}
One will obtain the following factorization $\bbA=\bbE\tilde\bbA$ where
 \ben
 \bbE &=& \left(\begin{array}{cccccc}
  \bbI_A & -\bbG^\top\bbL_\phi^{-1}  &0 &0  &0 &0 \\
  0 & \bbI_\phi &0 &0  &0 &0\\
  0 & 0 &\bbI_H & -\bbD^\top\bbL_r^{-1}  &0 &0\\
  0 & 0 &0 &\bbI_r  &0 &0\\
  0 & 0 &0 &0  &\bbI_u &0\\
  0 & 0 &0 &0  &0 &\bbI_p\\
 \end{array}\right),\\
 \tilde\bbA &=& \left(
\begin{array}{cccccc}
\tilde\bbC & \bbG^\top & \bbK  & 0         & 0         & 0 	\\
\bbG & 0         & 0     & 0         & 0         & 0 	\\
0    & 0         & \tilde\bbH  & \bbD^\top & \bbJ^\top & 0 	\\
0    & 0         & \bbD  & 0         & 0         & 0 	\\
0    & 0         & -\bbJ & 0         & \bbF      & \bbB^\top \\
0    & 0         & 0     & 0         & \bbB      & 0
\end{array}
\right).
 \een
The diagonal block matrices of $\bbE$ are all identity matrices of different sizes,
$\bbL_\phi$ is the stiffness matrix of $-\Delta$ on $Y_h$, $\bbL_r$ is the stiffness matrix of $-(R_m/\kappa)\Delta$ on $S_h$, and
\ben
\tilde\bbC=\bbC + \bbG^\top \bbL_\phi^{-1} \bbG,\qquad
\tilde\bbH =\bbH + \bbD^\top \bbL_r^{-1} \bbD.
\een
$\tilde\bbC$ and $\tilde\bbH$ are called mass augmentation of $\bbC$ and $\bbH$ because we can use mass matrices to approximate
$\bbG^\top \bbL_\phi^{-1} \bbG$ and $\bbD^\top \bbL_r^{-1} \bbD$ \cite{greif2007, phi16}.
We further consider the LU factorization $\tilde\bbA=\bbL\bbU$ where
 \ben
 \bbL &=&
 \left(\begin{array}{cccccc}
  \bbI_A &0  &0 &0  &0 &0 \\
  \bbG\tilde{\bbC}^{-1} & \bbI &0 &0  &0 &0\\
  0 & 0 &\bbI_\phi &0  &0 &0\\
  0 & 0 &\bbD\tilde{\bbH}^{-1} &\bbI_H  &0 &0\\
  0 & 0 & -\bbJ\tilde\bbH^{-1}
    & -\bbJ\tilde\bbH^{-1}\bbD^\top\bbS_r^{-1}  &\bbI_u &0\\
  0 & 0 &0 &0  &\bbB\tilde\bbF^{-1} &\bbI_p\\
\end{array}\right) ,\\
\bbU &=&\left(
\begin{array}{cccccc}
\tilde{\bbC} & \bbG^\top & \bbK  & 0         & 0         & 0 	\\
0 & -\bbS_\phi        & -\bbG \tilde{\bbC}^{-1}\bbK     & 0         & 0         & 0 	\\
0    & 0         & \tilde{\bbH}  & \bbD^\top & \bbJ^\top & 0 	\\
0    & 0         & 0  & -\bbS_r           & -\bbD \tilde{\bbH}^{-1} \bbJ^\top       & 0 	\\
0    & 0         & 0 & 0  & \tilde\bbF   & \bbB^\top \\
0    & 0         & 0     & 0         & 0      & -\bbS_p
\end{array}
\right).
 \een
The diagonal block matrices of $\bbU$ are given by
 \ben
 &&\bbS_\phi = \bbG \tilde{\bbC}^{-1} \bbG^\top,\qquad
 \bbS_r = \bbD \tilde{\bbH}^{-1} \bbD^\top, \qquad
 \bbS_p = \bbB\tilde\bbF^{-1}\bbB^\top,\\
 &&\tilde\bbF := \bbF+\bbJ\tilde\bbH^{-1}\bbJ^\top
      -(\bbD\tilde{\bbH}^{-1}\bbJ^\top)^\top \bbS_r^{-1}
      (\bbD\tilde{\bbH}^{-1}\bbJ^\top).
 \een
It is easy to see that $\bbA (\bbE\bbU)^{-1} = \bbE\bbL\bbE^{-1}$.
Since $\bbL$ only has unit eigenvalues, we call $(\bbE\bbU)^{-1}$ an {\em ideal preconditioner} of $\bbA$.

Direct calculations show
 \ben
 \bbE\bbU =\left(
\begin{array}{cccccc}
\tilde{\bbC} & \bbX_{12} & \bbX_{13} & 0  & 0 & 0 	\\
0 & -\bbS_\phi        & \bbX_{23}     & 0         & 0         & 0 	 \\
0    & 0   & \tilde{\bbH}  & \bbX_{34} & \bbX_{35}   & 0 	\\
0    & 0         & 0  & -\bbS_r  & \bbX_{45}   & 0 	\\
0    & 0         & 0 & 0  & \tilde\bbF   & \bbB^\top \\
0    & 0         & 0     & 0         & 0      & -\bbS_p
\end{array}
\right) ,
 \een
where
 \ben
 \begin{array}{lll}
 \bbX_{12}=\bbG^\top+\bbG^\top\bbL_\phi^{-1}\bbS_\phi,
 &\quad
    \bbX_{13}=\bbK-\bbG^\top\bbL_\phi^{-1}\bbX_{23},
 &\quad
    \bbX_{23} = -\bbG \tilde{\bbC}^{-1}\bbK,
 \vspace{1mm}\\
 \bbX_{34} =\bbD^\top+\bbD^\top\bbL_r^{-1}\bbS_r,
 &\quad
    \bbX_{35} = \bbJ^\top-\bbD^\top\bbL_r^{-1}\bbX_{45},
 &\quad
    \bbX_{45}=-\bbD \tilde{\bbH}^{-1} \bbJ^\top.
 \end{array}
 \een

\subsection{Practical preconditioner of $\bbA$}

The block entries of $\bbE\bbU$ is too complex to provide a practical preconditioner. We simplify them based on heuristic analysis.

Remember that the block matrices of $\bbA$ are algebraic representations of differential operators or multiplication operators appearing in the MHD system,
e.g.,
\begin{equation}\label{eq:tC}
 \tilde\bbC \Leftrightarrow \curl\curl +
    \nabla \left(-\Delta|_{Y_h}\right)^{-1}(-\Div) ,
    \quad  \bbG \Leftrightarrow -\Div \quad
    \hbox{on}\;\; \BD_h,
\end{equation}
where $(-\Div)$ is understood as the dual operator of $\nabla|_{Y_h}$. Similarly
\begin{equation}\label{eq:tH}
 \tilde\bbH \Leftrightarrow \kappa R_m^{-1}\left[\curl\curl +
    \nabla \left(-\Delta|_{S_h}\right)^{-1}(-\Div)\right] ,
    \quad  \bbD \Leftrightarrow -\Div \quad
    \hbox{on}\;\; \BW_h,
\end{equation}
where $(-\Div)$ is understood as the dual operator of $\nabla|_{S_h}$.
From \cite{greif2007}, we have the following spectral equivalences of matrices
 \begin{equation}\label{mat:CH}
 \tilde\bbC \sim\hat\bbC:=\bbC+\bbM_A,\quad
 \tilde\bbH \sim\hat\bbH:=\bbH+\kappa R_m^{-1}\bbM_H,\quad
 \bbS_\phi \sim\bbL_\phi,\quad \bbS_r \sim \bbL_r,
 \end{equation}
where $\bbM_A$, $\bbM_H$ are mass matrices on $\BD_h$ and $\BW_h$ respectively. This inspires us to make the replacements for off-diagonal blocks of $\bbE\bbU$
 \begin{equation}\label{mat:X12}
 \bbX_{12} \approx 2\bbG^\top,\qquad
 \bbX_{34} \approx 2\bbD^\top.
 \end{equation}

Now we shall follow the arguments of \cite{li17} to estimate $\bbX_{23}$ and $\bbX_{45}$.
Note that $\bbK$ and $\bbJ$ are algebraic representations of two multiplication operators
\begin{equation}\label{eq:KJ}
\bbK \Leftrightarrow (-\curl \BH_h^n)\quad
    \hbox{on}\;\;\BD_h,\qquad
\bbJ^\top \Leftrightarrow
    \kappa\curl(\curl\BA_h^{n-1}\times \Bu_h^n)
    \quad \hbox{on}\;\;\BW_h.
\end{equation}
Formally we have the identity $(-\Div)(\curl\curl+\nabla\Delta^{-1}\Div)=-\Div$ on $\BC_0^\infty(\Omega)$. This shows that
$(-\Div)(\curl\curl+\nabla\Delta^{-1}\Div)^{-1}
=-\Div$ on $\BC_0^\infty$. Now from \eqref{eq:tC} and \eqref{eq:KJ}, we deduce heuristically that
\begin{equation}\label{mat:X23}
\bbX_{23} = -\bbG \tilde{\bbC}^{-1}\bbK
\approx -\bbG\bbK \approx 0.
\end{equation}
Similarly, from \eqref{eq:tC} and \eqref{eq:KJ}, we deduce heuristically that
\begin{equation}\label{mat:X45}
\bbX_{45} = -\bbD \tilde{\bbH}^{-1} \bbJ^\top
 \approx -\bbD\bbJ^\top \approx 0.
\end{equation}
Based on \eqref{mat:X23} and \eqref{mat:X45}, we also have the approximations
\begin{equation}\label{mat:X13}
 \bbX_{13} \approx \bbK,\qquad
 \bbX_{35} \approx \bbJ^\top.
\end{equation}

Now it is left to consider the Navier-Stokes block of $\bbE\bbU$.
From \eqref{mat:CH} and \eqref{mat:X45}, we have the following approximations
\ben
\tilde\bbF \approx \bbF + \bbJ\tilde\bbH^{-1}\bbJ^\top
\approx  \bbF + \bbJ\hat\bbH^{-1}\bbJ^\top,\quad
\bbS_p =\bbB\tilde\bbF^{-1}\bbB^\top
\approx \bbB\left(\bbF + \bbJ\hat\bbH^{-1}\bbJ^\top\right)^{-1}\bbB^\top.
\een
Here $\bbJ\hat\bbH^{-1}\bbJ^\top$ stands for coupling term between the fluid and electromagnetic field.
In \cite{li17}, Li and Zheng derived approximate Schur complements for $\bbF + \bbJ\hat\bbH^{-1}\bbJ^\top$ and $\bbS_p$ in the case that the Navier-Stokes equations are solved by $\BP_2$--$P_1$ Taylor-Hood finite elements.
Here we apply their results directly to our case of $(\Bu_h^n,p_h^n)\in\BV_h\times Q_h$. In \cite[Section 3.3]{li17},
Li and Zheng suggested to approximate $\bbF + \bbJ\hat\bbH^{-1}\bbJ^\top$ and $\bbS_p$ as follows
 \begin{equation}\label{mat:F}
 \tilde\bbF \approx
 \bbF+\bbJ\hat\bbH^{-1}\bbJ^\top\approx \bbS_u, \qquad
 \bbS_p \approx \bbB\bbS_u^{-1}\bbB^\top
 \approx (R_e^{-1}+\alpha)^{-1}\bbM_p,
 \end{equation}
where $\bbS_u$ is the stiffness matrix associated with the bilinear form
\ben
\mathscr{A}_u(\Bw,\Bv):=
\mathscr{A}_1(\Bw,\Bv)+\mathscr{O}_h(\Bu_h^{n-1};\Bw,\Bv)+
\kappa{R}_m\SP{\curl\BA_h^{n-1}\times\Bw}{\curl\BA_h^{n-1}\times\Bv},
\een
which means that we use $\kappa{R}_m\SP{\curl\BA_h^{n-1}\times\Bw}{\curl\BA_h^{n-1}\times\Bv}$ to
approximate $\bbJ\hat\bbH^{-1}\bbJ^\top$ \cite{li17}.
Here $\bbM_p$ is the mass Matrix on finite element space $Q_h$, and $\alpha$ is the grad-div stabilization parameter.
We refer to \cite{ben06, benzi2011FOV, ben11, farrell2019} for more details about augmented Lagrangian preconditioners
for solving Navier-Stokes equations.

Finally, using the approximations \eqref{mat:CH}--\eqref{mat:F} in $\bbE\bbU$, a practical preconditioner of $\bbA$ can be defined by the inverse of
\begin{align}\label{pre}
\bbP = \left(
\begin{array}{cccccc}
\hat{\bbC} & 2\bbG^\top & \bbK & 0  & 0 & 0 	\\
0 & -\bbL_\phi & 0     & 0         & 0         & 0 	 \\
0    & 0   & \hat{\bbH}  & 2\bbD^\top & \bbJ^\top & 0 	\\
0    & 0   & 0  & -\bbL_r  & 0   & 0 	\\
0    & 0   & 0  & 0  & \bbS_u   & \bbB^\top \\
0    & 0   & 0  & 0  & 0    & -(R_e^{-1}+\alpha)^{-1}\bbM_p
\end{array}
\right).
\end{align}

\subsection{A preconditioned GMRES algorithm}
Based on $\bbP$, we propose a preconditioned GMRES method for solving \eqref{eq:algeb}.
In each GMRES iteration, one needs to solve the system of algebraic equations
\begin{align}\label{eq:pre}
\bbP\Ve=\Vr,
\end{align}
where $\Ve$ is the correction vector and
$\Vr$ the residual vector calculated from last iteration.
Now we present the algorithm for solving an approximate solution of (\ref{eq:pre}). \vspace{2mm}

\begin{algorithm}\label{alg-pre}
{\sf
Set the tolerance $\varepsilon_0=10^{-3}$ and write
 \ben
 \Ve=(\Ve_A,\Ve_\phi,\Ve_H,\Ve_r,\Ve_u,\Ve_p)^\top,\qquad \Vr=(\Vr_A,\Vr_\phi,\Vr_H,\Vr_r,\Vr_u,\Vr_p)^\top.
 \een
The approximate solution of \eqref{eq:pre} is computed in six steps below. In each step, the algebraic problem is solved iteratively until the relative residual is less than $\varepsilon_0$.
\vspace{1mm}

\begin{enumerate}[leftmargin=6mm]
\item[1.] Solve $\bbM_p\Ve_p = - (R_e^{-1}+\alpha)\Vr_p$ by the CG method with diagonal preconditioner.  \vspace{1mm}

\item[2.] Solve $\bbS_u\Ve_u=\Vr_u-\bbB^\top\Ve_p$ by the GMRES method with additive Schwarz preconditioner (cf. \cite{cai99}).   \vspace{1mm}

\item[3.] Solve $\bbL_r\Ve_r = -\Vr_r$ by the CG method with algebraic multigrid solver (cf. \cite{hen02}).   \vspace{1mm}

\item[4.] Solve $\hat\bbH\Ve_H=\Vr_H-2\bbD^\top\Ve_r-\bbJ^\top\Ve_u$ by the CG method with auxiliary space preconditioner (cf. \cite{hip07}). \vspace{1mm}

\item[5.] Solve $\bbL_\phi\Ve_\phi = -\Vr_\phi$ by the CG method with algebraic multigrid solver.  \vspace{1mm}

\item[6.] Solve $\hat{\bbC}\Ve_A=\Vr_A-2\bbG^\top\Ve_\phi-\bbK\Ve_H$ by the CG method with auxiliary space preconditioner.
\end{enumerate}}
\end{algorithm}
\vspace{2mm}

We remark that the total number of iterations for solving $\bbA\Vx = \Vb$ is insensitive to the choice of tolerance $\varepsilon_0$ in Algorithm~\ref{alg-pre} when $\varepsilon_0\le 10^{-2}$.
We should confess that the additive Schwarz preconditioner adopted in Step 2 of Algorithm~\ref{alg-pre} is not optimal.
More efficient preconditioners for $\bbS_u$ are important to improve the overall efficiency and will be our future work.
We refer the readers to recent work \cite{farrell2019} on three-dimensional stationary Navier-Stokes equations using augmented Lagrangian block preconitioner, where an optimal geometrical multigrid method is developed.
However an extension of the multigrid techniques for the solving
of the Step 2 problem, where $\BH(\Div)$-conforming element
is used for the velocity field $\Bu_h$, is not a trivial thing.

\section{Numerical experiments}
In this section, we report several numerical experiments to show convergence orders of discrete solutions and
to demonstrate the performance of the preconditioned GMRES solver.
The finite element method and the discrete solver are implemented
on the finite element package "Parallel Hierarchical Grid" (PHG) \cite{zhang09}.

For solving the nonlinear problem \eqref{weakh}, the Picard iterations stop whenever the criterion is reached
\ben
\Theta(\Bu_h^n) +\Theta(\BH_h^n)+\Theta(\BA_h^n)< \delta,
\een
where $\Theta(\Bw_h^n)=\NLtwov{\Bw_h^n-\Bw_h^{n-1}}\NLtwov{\Bw_h^n}^{-1}$
for $\Bw_h^n=\Bu_h^n$, $\BH_h^n$, and $\BA_h^n$. Here $\delta$ is the tolerance for Picard's iterations.
For solving the linear problem \eqref{eq:algeb}, let $\Vx^{(k)}$, $k\ge 0$, be the approximate solution at $k^{\text{th}}$ GMRES iteration and let $\Vr^{(k)}=\Vb-\bbA\Vx^{(k)}$ be the residual.
The iterations stop whenever the criterion is reached
\begin{align*}
\big\|\Vr^{(k)}\big\|_2\le \varepsilon
\big\|\Vr^{(0)}\big\|_2,
\end{align*}
where $\varepsilon$ is the tolerance for the GMRES solver.
The maximal iteration number for the GMRES solver is set by $200$ without restart
and right preconditioning algorithm is adopted here.

Throughout this section, we set the penalty parameter in \eqref{eq:Ah} by $\gamma = 10$
and the grad-div parameter in $\mathscr{A}_1$ by $\alpha=1$, except
for Example~\ref{ex:alpha} where the sensitivity of the solver to $\alpha$ is tested. The domain is chosen as $\Omega = [0,1]^3$.
We choose 5 quasi-uniform meshes of $\Omega$ by successive refinements. The information of the meshes is listed in
Table~\ref{tab1}.

\begin{table}[htbp!]
  \captionsetup{position=top}
  \caption{Five successively refined meshes.}  \label{tab1}
  \centering
  \begin{tabular}{|c|c|r|r|r|}
  \hline
  Mesh & $h$ & DOFs for $(\BA_h,\psi_h)$ & DOFs for $(\BH_h,r_h)$ & DOFs for $(\Bu_h,p_h)$ \\ \hline
  $\Ct_1$	 & 1.732  & 65     & 65      & 60 \\ \hline
  $\Ct_2$	 & 0.866  & 321    & 321     & 408 \\ \hline
  $\Ct_3$	 & 0.433  & 1,937  & 1,937   & 2,976\\ \hline
  $\Ct_4$	 & 0.217  & 13,281 & 13,281  & 22,656\\ \hline
  $\Ct_5$	 & 0.108  & 97,985 & 97,985  & 176,640 \\ \hline
\end{tabular}
\end{table}

\begin{example}\label{ex:conv}
This example is to investigate convergence orders of finite element solutions. The physical parameters are set by $R_e=R_m=\kappa=1$. The tolerances are set by $\delta=10^{-5}$ and $\varepsilon=10^{-6}$. The right-hand sides and the Dirichlet boundary conditions are chosen so that the true solutions are given by
\begin{align*}
&\BA=(\sin z,0,0)^\top,\quad \BH=(0,\cos z,0)^\top, \quad\psi={r}=0,\\
&\Bu=(\cos z, \sin(x+z),0)^\top,\qquad {p}=x+y-1.
\end{align*}
\end{example}

From Table \ref{tab2}-\ref{tab3}, we find that optimal convergence orders are obtained for physical quantities, namely, $\BA_h$, $\BH_h$, $\Bu_h$, and $p_h$, under their energy norms
\begin{align*}
\begin{array}{ll}
\NHcurl{\BA-\BA_h}\sim{O}(h),
&\quad\NHcurl{\BH-\BH_h}\sim{O}(h),\vspace{1mm}\\
\Ndg{\Bu-\Bu_h}\sim{O}(h),
&\quad \NLtwo{p-p_h}\sim{O}(h).
\end{array}
\end{align*}
Moreover, we also find that $\NLtwo{\Div\Bu_h}$ is negligible, compared with approximation errors.
The reason for $\Div\Bu_h\ne 0$ is due to the error from solving the system of linear algebraic equations \eqref{eq:algeb}, namely, the tolerance $\varepsilon =10^{-6}$.

\begin{table}[htbp!]
  \captionsetup{position=top} 
  \caption{Convergence orders of $\Bu_h$ and $p_h$.
  {\rm (Example~\ref{ex:conv})}}
  \label{tab2}
  \centering
  \begin{tabular}{|c|c|c|c|c|c|}
  \hline
  Mesh  & $\Ndg{\Bu-\Bu_h}$ & Order & $\NLtwo{p-p_h}$ & Order & $\NLtwo{\Div\Bu_h}$  \\ \hline
  $\Ct_1$& 6.376e-01 & ---    & 1.227e+01 & ---    & 3.701e-10\\ \hline
  $\Ct_2$& 2.411e-01 & 1.403  & 4.040e-01 & 4.925  & 2.149e-09\\ \hline
  $\Ct_3$& 1.203e-01 & 1.003  & 1.399e-01 & 1.530  & 5.067e-09\\ \hline
  $\Ct_4$& 5.865e-02 & 1.036  & 5.221e-02 & 1.422  & 2.123e-08\\ \hline
  $\Ct_5$& 2.874e-02 & 1.029  & 2.215e-02 & 1.237  & 6.088e-08\\ \hline
 \end{tabular}
\end{table}

\begin{table}[htbp!]
  \captionsetup{position=top} 
  \caption{Convergence orders of $\BA_h$ and $\BH_h$.
  {\rm (Example~\ref{ex:conv})}}
  \label{tab3}
  \centering
  \begin{tabular}{|c|c|c|c|c|}
  \hline
  Mesh  & $\NHcurl{\BA-\BA_h}$ & Order & $\NHcurl{\BH-\BH_h}$ & Order  \\ \hline
  $\Ct_1$& 1.190e-01 & ---    & 2.114e-01 & ---   \\ \hline
  $\Ct_2$& 5.937e-02 & 1.003  & 9.883e-02 & 1.097 \\ \hline
  $\Ct_3$& 2.942e-02 & 1.013  & 4.862e-02 & 1.023 \\ \hline
  $\Ct_4$& 1.457e-02 & 1.014  & 2.410e-02 & 1.013 \\ \hline
  $\Ct_5$& 7.234e-03 & 1.010  & 1.200e-02 & 1.001 \\ \hline
  \end{tabular}
\end{table}
\vspace{5mm}

\begin{example}[Driven cavity flow]\label{ex:cavity}
The purpose of this example is to demonstrates the optimality of the preconditioned GMRES method for solving a benchmark problem. The external force in the momentum equation is set by $\Bf=\mathbf{0}$. The boundary conditions are set by
\begin{align*}
\BA=(0,0,-y),\quad\BH=(-1,0,0),\quad\Bu=(v,0,0)\quad
\text{on}\;\;\partial\Omega,
\end{align*}
where $v\in C[0,1]$ and satisfies
\[
v(x,y,1)=1,\qquad
v(x,y,z)=0\quad\forall\,z\in[0,1-h].
\]
\end{example}

We fix $\kappa=1$, $R_m=10$ and demonstrate the optimality of the solver for $R_e=1,10,100$.
The tolerances are set by $\delta=10^{-4}$ and $\varepsilon=10^{-5}$.
Let $N_{\text{picard}}$ denote the number of Picard iterations and let $N_{\text{gmres}}$ denote the average number of preconditioned GMRES iterations for solving \eqref{eq:algeb}.
From Table~\ref{tab4}, we find that the number of GMRES iterations is quasi-uniform to the meshes for each fixed $R_e$. Moreover,  the preconditioned GMRES method is robust to $R_e$.

\begin{table}[htbp!]
\captionsetup{position=top} 
\caption{Robustness and quasi-optimality of the preconditioned GMRES solver. {\rm (Example~\ref{ex:cavity})}}
\label{tab4}
\centering
\begin{tabular}{|c|c|c|c|}
\hline
  & \multicolumn{3}{c|}{$N_{\text{gmres}}(N_{\text{picard}})$} \\ \hline
 \diagbox{Mesh}{$R_e$}& {$1$} &  {$10$} & {$100$} \\ \hline
$\Ct_1$& \;\;9 (5)   & 9 (4)  & \;\;7 (5)  \\ \hline
$\Ct_2$& 12 (5)  & 8 (4)  & \;\;8 (8)  \\ \hline
$\Ct_3$& 10 (5)  & 8 (5)  & 10 (7)  \\ \hline
$\Ct_4$& 11 (6)  & 8 (5)  & 14 (7)  \\ \hline
$\Ct_5$& \;\;9 (6)   & 9 (6)  & 15 (7)  \\ \hline
\end{tabular}
\end{table}
\vspace{5mm}

\begin{example}[Robustness]\label{ex:robust}
This example investigates the robustness of the solver to $R_e$ and $R_m$ by the driven cavity flow in Example~\ref{ex:cavity}.
\end{example}

\begin{figure}[htbp!]
  \centering
  \includegraphics[width=0.32\textwidth]{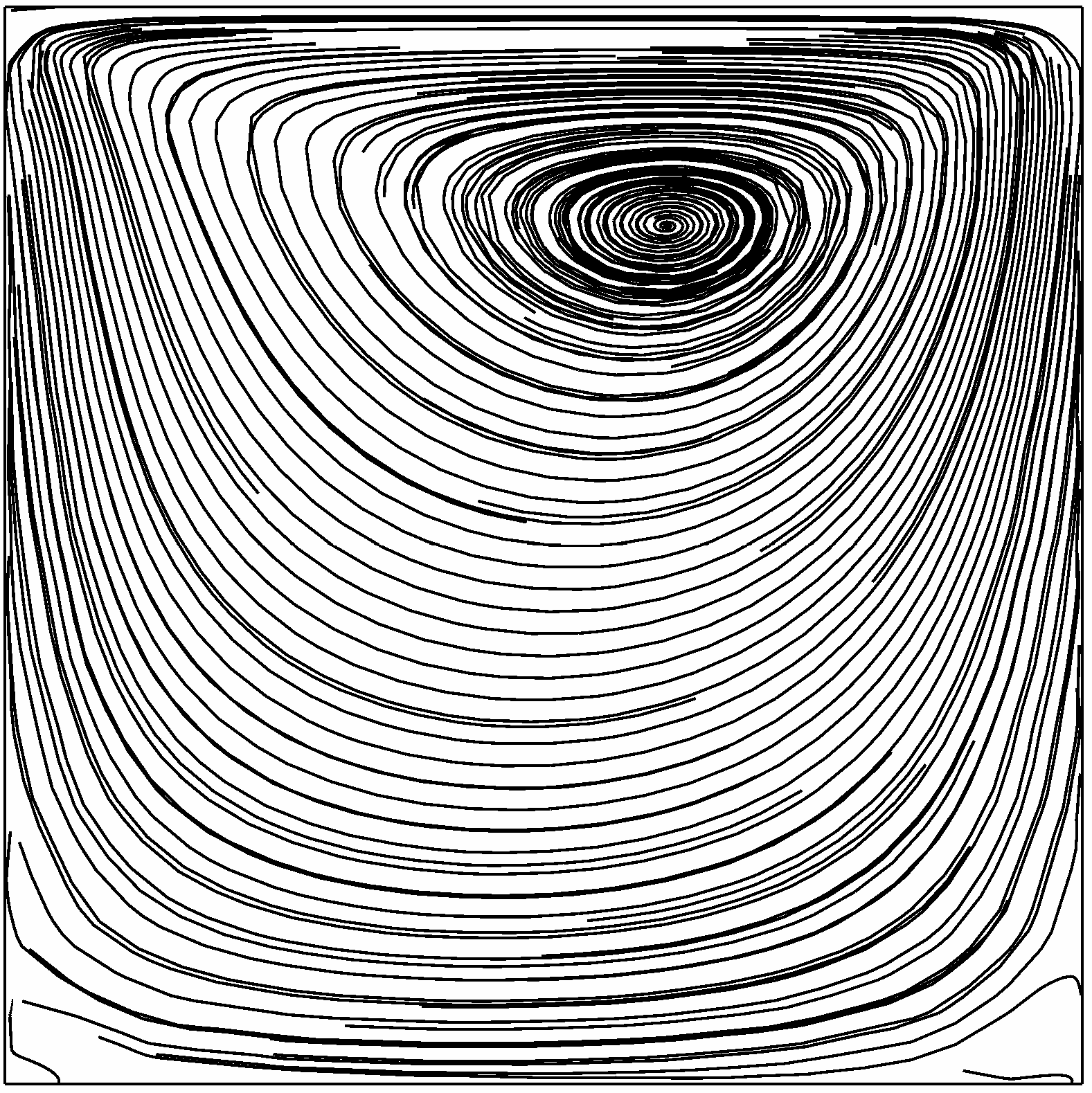}
  \includegraphics[width=0.32\textwidth]{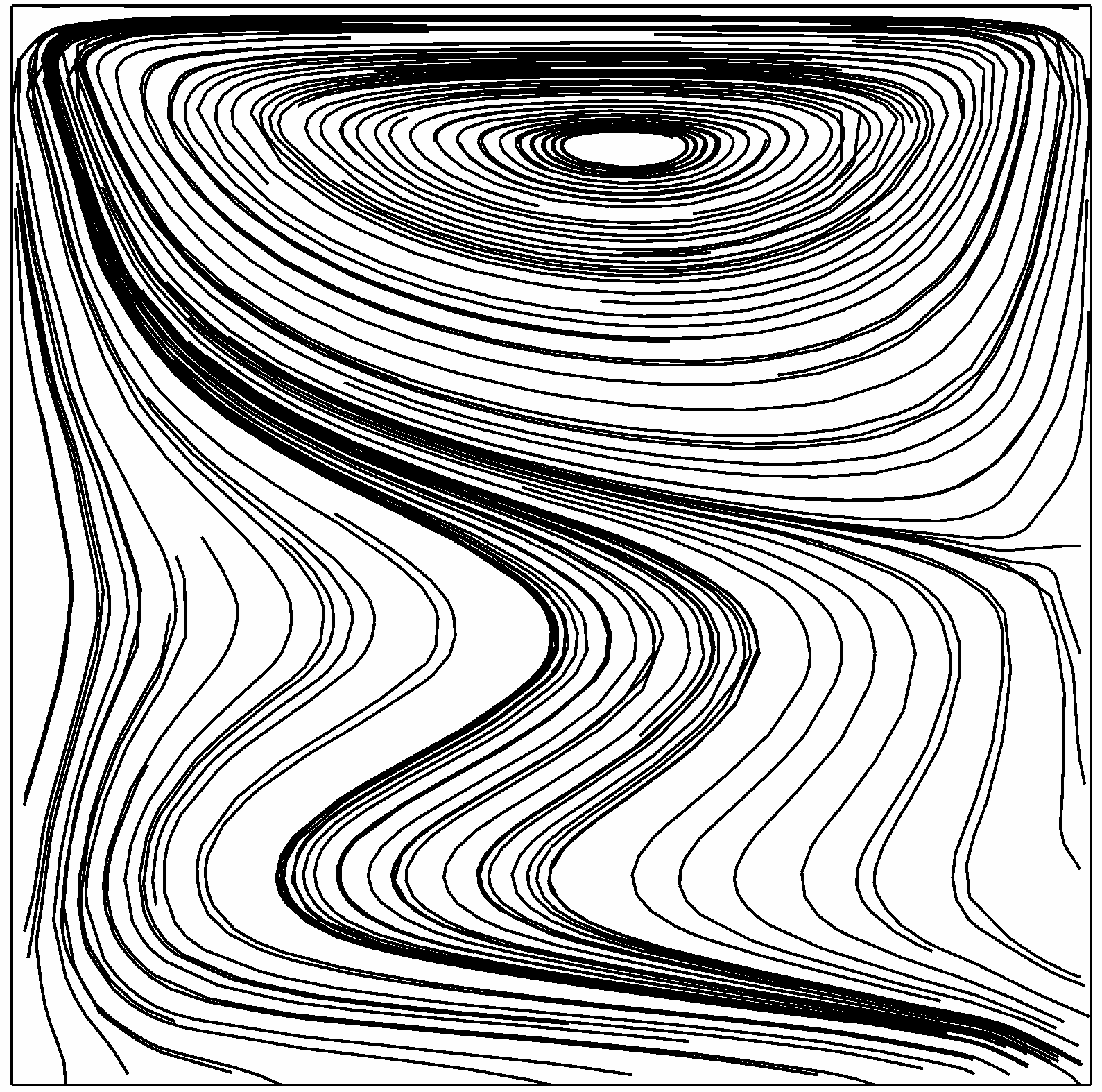}
  \includegraphics[width=0.32\textwidth]{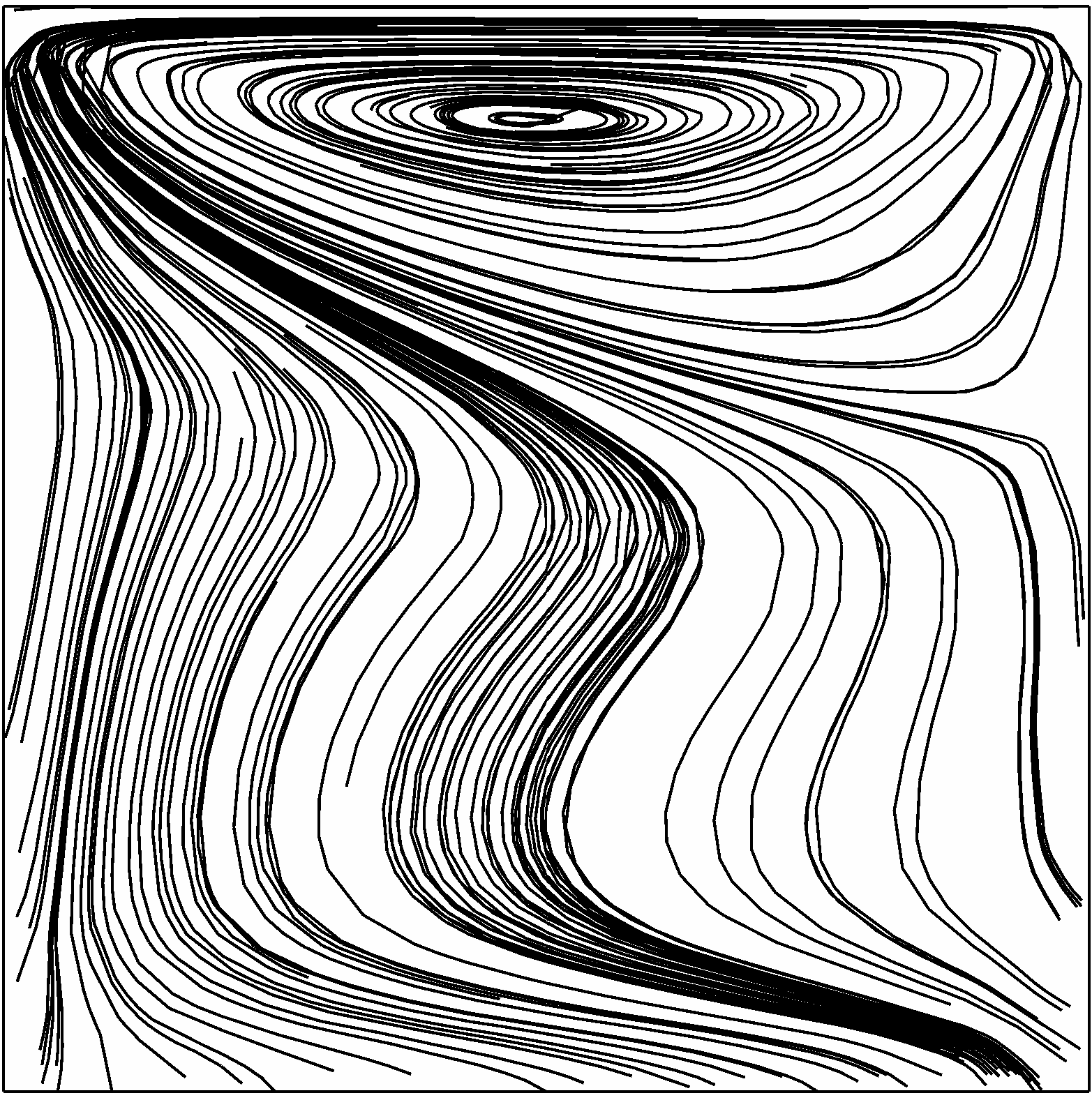}
  \caption{Projections of the streamlines of $\Bu_h$ on the cross
  section $y=0.5$ with $R_e=100$ and $\kappa=1$ $($from left to the right $R_m=1,10,50)$.
  {\rm(Example~\ref{ex:robust})}}\label{fig3}
\end{figure}

We choose $\Ct_5$ as the computational mesh and set $\kappa=\alpha=1$.
The tolerances are $\delta=10^{-4}$ and $\varepsilon=10^{-5}$.
From Table~\ref{tab5}, we find that, for small $R_m$, the solver for linear system is robust with respect to $R_e$,
while for large $R_m$, the number of GMRES iterations grows slightly.
Another observation is that with large $R_m$,
Picard's method is inefficient for solving the nonlinear problem \eqref{weakh}.
Unfortunately, the present discretization using upwinding in the convection term $\Bw\cdot\nabla\Bu$,
which makes the Newton's method difficult to use.
In the future, acceleration techniques in optimization field for nonlinear iteration
can be incorporated to improve the nonlinear convergence rate.
In Figure~\ref{fig3}, we depict the streamlines of $\Bu_h$ projected onto the cross-section at $y = 0.5$ for different values of $R_m$.

\begin{table}[htbp!]
  \captionsetup{position=top} 
  \caption{Sensitivity to $R_e$ and $R_m$.
  {\rm(Example~\ref{ex:robust})}}\label{tab5}
  \centering
  \begin{tabular}{|c|c|c|c|c|}
  \hline
  & \multicolumn{4}{c|}{$N_{\text{gmres}}(N_{\text{picard}})$}\\
  \hline
  \diagbox{$R_e$}{$R_m$} & 1 & 20 & 40 & 60        \\ \hline
  1 	& 6 (4)  & 10 (11)    & 12 (26) & 14 ($>$100) \\ \hline
  10	& 5 (5)  & 11 (10)    & 13 (30) & 17 ($>$100) \\ \hline
  100	& \;\;8 (10) & 20 (8)\;\; & 31 (44) & 39 ($>$100) \\ \hline
  \end{tabular}
\end{table}

Next we choose $\Ct_5$ as the computational mesh and set $\alpha=1,R_e=100$.
We examine the effect of the magneto-fluid coupling term on the performance of the preconditioner.
Remember from \eqref{mat:F} that the approximate matrix of $\tilde\bbF$ is $\bbS_u\approx\bbF+\bbJ\hat{\bbH}^{-1}\bbJ^\top$.
The data in Table~\ref{tab6} is the result with the proposed preconditioner,
while in Table~\ref{tab7} we ignore the coupling term $\bbJ\hat{\bbH}^{-1}\bbJ^\top$ and only use $\bbF$ as the approximate matrix.
From Table~\ref{tab7}, we see that the numbers of the GMRES increase considerably for large $\kappa$ or $R_m$.
\begin{table}[htbp!]
  \captionsetup{position=top} 
  \caption{Performance for different $\kappa$ and $R_m$ with approximation $\bbS_u$.
  {\rm(Example~\ref{ex:robust})}}\label{tab6}
  \centering
  \begin{tabular}{|c|c|c|c|}\hline
  &\multicolumn{3}{c|}{$N_{\text{gmres}}(N_{\text{picard}})$}\\\hline
  \diagbox{$\kappa$}{$R_m$}       &1        &20        &40                                  \\\hline
  \;1    &\multicolumn{1}{r|}{8\;(10)}      &20 (8)    &\multicolumn{1}{r|}{\;\;\;31\;(44)} \\\hline
  20     &\multicolumn{1}{r|}{24\;(7)\,\;}  &70 (6)    &\multicolumn{1}{r|}{\,\,90\;(7)}    \\\hline
  40     &\multicolumn{1}{r|}{32\;(7)\,\;}  &87 (6)    &\multicolumn{1}{r|}{109 (6)}        \\\hline
  60     &\multicolumn{1}{r|}{38\;(7)\,\;}  &99 (6)    &\multicolumn{1}{r|}{125 (6)}        \\\hline
  \end{tabular}
\end{table}

\begin{table}[htbp!]
  \captionsetup{position=top} 
  \caption{Performance for different $\kappa$ and $R_m$ with approximation $\bbF$.
  {\rm(Example~\ref{ex:robust})}}\label{tab7}
  \centering
  \begin{tabular}{|c|c|c|c|}  \hline
  &\multicolumn{3}{c|}{$N_{\text{gmres}}(N_{\text{picard}})$} \\ \hline
  {\diagbox{$\kappa$}{$R_m$}}   &1         &20            &40                                      \\\hline
  1     &\multicolumn{1}{r|}{7\;(10)}      &\,\,26\;\,(7) &\multicolumn{1}{r|}{37  (47)\;\;\;\,}           \\\hline
  20 	&\multicolumn{1}{r|}{24\;(7)\,\;}  &119 (6)       &\multicolumn{1}{r|}{171 (8)\quad\;\;\,} \\\hline
  40	&\multicolumn{1}{r|}{35\;(7)\,\;}  &161 (6)       &\multicolumn{1}{r|}{171 (7)\quad\;\;\,} \\\hline
  60	&\multicolumn{1}{r|}{43\;(7)\,\;}  &167 (6)       &\multicolumn{1}{r|}{175 (8)\quad\;\;\,} \\\hline
  \end{tabular}
\end{table}
\vspace{5mm}

\begin{example}\label{ex:alpha}
This example investigates the sensitivity of the preconditioner to the grad-div stabilization parameter $\alpha$ by the driven cavity flow in Example~\ref{ex:cavity}.
\end{example}

We fix $R_m=1$ and $\kappa=100$ and investigate the performance of the GMRES solver to $R_e$ and $\alpha$. From Table~\ref{tab8}, we find that,
\begin{itemize}[leftmargin=6mm]
  \item for $\alpha\ge 0.5$, the convergence of the solver is not sensitive to $\alpha$,
  \item while for $\alpha=0$, the number of GMRES iterations increases fast with $R_e$.
\end{itemize}
We conclude that the grad-div stabilization plays an important role in the performance of the block preconditioner.
\begin{table}[htbp!]
  \captionsetup{position=top} 
  \caption{Sensitivity to $R_e$ and $\alpha$.
  {\rm(Example~\ref{ex:alpha})}}\label{tab8}
  \centering
  \begin{tabular}{|c|c|c|c|c|c|}  \hline
  &\multicolumn{5}{c|}{$N_{\text{gmres}}(N_{\text{picard}})$} \\\hline
  \diagbox{$R_e$}{$\alpha$} & 0 & 0.5 & 1 & 10  & 100      \\\hline
  1 & \;\;13 (4) & 11 (4)  & 10 (4) & \;\;9 (4)&\;\;7 (5)  \\\hline
  10& \;\;29 (5) & 25 (5)  & 21 (5) & 15 (4) &  12 (4)     \\\hline
  100&   105 (6) & 46 (6)  & 43 (6) & 34 (6) &  29 (6)     \\\hline
  \end{tabular}
\end{table}


\section{Conclusions}
In this paper, we propose a monolithic constrained transport finite element method for stationary incompressible MHD equations.
The discrete velocity, discrete current density, and discrete magnetic induction are all divergence-free in the momentum equation,
especially in the \textit{Lorentz force}.
Based on an augmented Lagrangian block preconditioner,
we also develop a preconditioned GMRES solver for the linearized system of algebraic equations in every Picard iteration.
Although the present work only presents first-order discretization for $\Bu$, $\BH$, and $\BA$,
the method can be extended to high-order finite elements straightforwardly.
The monolithic manner can be applied to time-dependent MHD equations to develop fully implicit method,
which will permit large time-step length and stable long time simulation compared with explicit method.

\section*{Acknowledgments}
The computations were (partly) done on the high performance computers of State Key Laboratory of Scientific and Engineering Computing,
Chinese Academy of Science.

Lingxiao Li was supported by National Natural Science Foundation of China under Grant 11901042.
Weiying Zheng was supported in part by the National Science Fund for Distinguished
Young Scholars 11725106 and by China NSF grant 11831016.

\end{document}